\pdfoutput=1
\documentclass{article}

\PassOptionsToPackage{numbers,square,sort&compress}{natbib}
\usepackage[margin=1in]{geometry}
\usepackage{natbib}
\usepackage[utf8]{inputenc}
\usepackage[T1]{fontenc}
\usepackage{url}
\usepackage{booktabs,tabularx,array}
\usepackage{amsmath, amsfonts, amssymb}
\usepackage{float}
\usepackage{graphicx}
\usepackage{microtype}
\usepackage[section]{placeins}
\usepackage{xcolor}
\usepackage{hyperref}
\date{}

\title{UNIDE: Cross-PDE Learning for Noise-Robust Discovery of Governing Equations}

\author{
Xiaoyu Xie\\
Independent Researcher\\
\texttt{xiaoyuxie.vico@gmail.com}}

\begin{document}

\maketitle

\begin{abstract}
Governing partial differential equations (PDEs) provide compact and interpretable descriptions of physical systems. Classical methods, such as sparse and symbolic regression, fit a new equation to each observed field: what is learned on one system is not carried to the next, and observation noise corrupts the derivative features they rely on. This paper introduces UNIDE (Unified Identification of Differential Equations), a data-driven approach that learns across PDE families to identify governing equations from new noisy trajectories. It integrates three components: a self-supervised denoiser trained without clean targets, an encoder shared across PDE families that predicts which candidate terms are active, and weak-form least squares with residual-based pruning that estimates their coefficients. We demonstrate UNIDE on synthetic one-dimensional trajectories from six PDE families with known term structures, including the Burgers and Kuramoto--Sivashinsky equations. We train one model per noise level and systematically compare UNIDE with SINDy, WSINDy, and PySR. The results show that UNIDE recovers the exact set of governing terms more often than all three at every noise level from $1\%$ to $80\%$. At $80\%$ noise, where the noise standard deviation is $0.8$ times that of the field, it recovers the exact terms for $92\%$ of held-out trajectories, whereas the best baseline recovers about $1\%$. UNIDE is also efficient: because its neural components are reused rather than refit, it identifies each new trajectory in $0.45$ seconds, $2.7$ times faster than WSINDy and $51$ times faster than PySR. These results demonstrate the potential of learning across PDE families for efficient equation identification under substantial observation noise.
\end{abstract}

\section{Introduction}

\begin{figure}[!htbp]
\centering
\includegraphics[width=\linewidth]{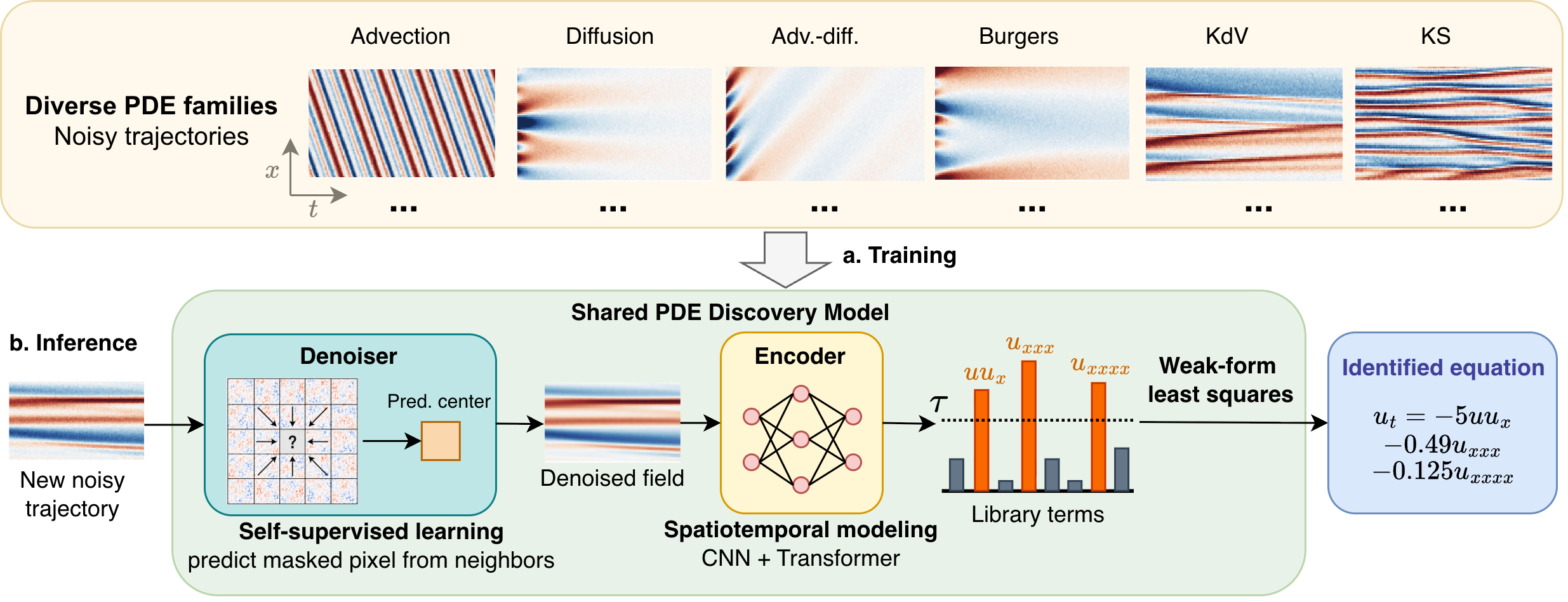}
\caption{The UNIDE method. The upper row illustrates training across six PDE families; the denoiser and encoder are trained sequentially, with a separate model for each observation-noise level. The lower row illustrates inference: the denoiser maps a noisy trajectory to a denoised field, the encoder predicts candidate-term probabilities, and weak-form least squares with residual pruning estimates the equation. The field images, probability bars, and displayed equation are schematic, not measured results.}
\label{fig:overview}
\end{figure}

Identifying the governing partial differential equation (PDE) of an observed system provides a compact and interpretable description of its dynamics: each term represents a distinct physical process, and the equation supports prediction beyond the observed regime. Despite the diversity of physical systems, most governing PDEs are built from a small set of recurring mechanisms, such as advection, diffusion, dispersion, and reaction, each of which leaves a characteristic signature in space--time data. For example, advection translates a profile without changing its shape, diffusion smooths it, and dispersion spreads it into trains of waves with different wavelengths. This recurrence suggests that equation discovery can be treated as a transferable inference problem rather than solved independently for each system: a model trained across many systems could learn these signatures and recover the governing equation directly from a new observation.

Identifying governing equations from noisy measurements nevertheless remains challenging. Established solutions are sparse regression and symbolic regression. SINDy and PDE-FIND~\citep{brunton2016discovering,rudy2017data} fit sparse combinations of candidate terms, and related approaches incorporate dimensional structure~\citep{xie2022data}; PySR~\citep{cranmer2023pysr} searches symbolic expressions with genetic programming. However, these methods have two well-known limitations. First, they rely on derivative features estimated from the data, and differentiation amplifies observation noise and impairs the selection of terms; weak-form~\citep{messenger2021weak} and ensemble~\citep{fasel2022ensemble} formulations reduce this sensitivity but do not remove it. Second, these methods optimize a new model for each observed space--time field, which we call a trajectory. Because each fit starts from scratch, experience gained on earlier systems cannot inform the next, and robustness to noise must be achieved one trajectory at a time.

These limitations could be overcome by learning across systems. The end-to-end transformer of \citet{kamienny2022end} predicts a full symbolic expression in one pass, but it is trained on algebraic functions of at most ten variables. ODEFormer~\citep{dascoli2024odeformer} extends this approach to ODEs, including noisy observations, but not to PDEs. PROSE-PDE~\citep{sun2024prosepde} couples multi-operator learning with symbolic PDE output, but its noisy-data experiments provide a symbolic expression or skeleton as an additional input alongside the noisy observations. FoundPDE~\citep{liu2026foundpde} pretrains a transformer on 20 one-dimensional PDE families, but its accuracy in recognizing the equation type falls from $97.5\%$ on clean data to $52.1\%$ under $10\%$ multiplicative noise.

Other approaches add learned priors or robustness to noise, but still run an optimization for each new system. EqGPT~\citep{xu2025eqgpt} pretrains an equation generator on handbook PDEs and recovers equations from noisy data, but for each new dataset it trains a neural surrogate of the observations and fine-tunes the generator in a generation--evaluation--optimization loop. Large language models supply similar priors~\citep{shojaee2025llmsr}, again with a search run separately for each dataset. DeepMoD~\citep{both2021deepmod} and Weak-PDE-Net~\citep{weakpdenet2026} target noisy data but optimize a new model for each observed system. On the MDBench benchmark~\citep{bideh2026mdbench}, DeepMoD recovers the full equation for 2 of 14 clean-data PDEs; it fails on four and does not complete the other eight within the benchmark's runtime, dimensionality, and sampling limits. Learning across systems and robustness to noise have thus been pursued largely separately. We ask whether one model can do both: learn the signatures of physical mechanisms across PDE families and identify the governing terms of a new trajectory at noise levels where per-trajectory methods fail.

In this study, we propose a unified data-driven approach, called \textbf{UNIDE} (Unified Identification of Differential Equations), that learns across PDE families to identify governing equations from noisy trajectories. The approach consists of three components (Figure~\ref{fig:overview}). First, a denoiser trained with a blind-spot objective~\citep{krull2019noise2void} maps the noisy trajectory to a denoised field without clean targets. Second, an encoder shared across PDE families predicts which candidate terms are active. Third, weak-form least squares estimates the coefficients of the selected terms, and greedy residual pruning removes selected terms whose deletion changes the fit little. The neural components are trained once per noise level and reused, without retraining, for every new trajectory. We demonstrate UNIDE in two settings: 900 held-out trajectories of the training systems at ten noise levels, and 360 trajectories of held-out systems with larger coefficients but the same term structures. We further compare it with sparse regression given the same denoiser, candidate terms, and pruning, which separates the contribution of learned term selection from that of denoising. At $80\%$ noise, UNIDE recovers the exact set of governing terms for $92\%$ of the held-out trajectories, whereas the best of the three baselines recovers about $1\%$, and the matched comparison attributes this advantage to learned term selection rather than to denoising alone.

\section{Problem Setting}
\label{sec:problem}

We study trajectories generated from one-dimensional periodic PDEs of the form
\[
u_t = \sum_{j=1}^{M} \xi_j \Theta_j(u), \qquad M=12,
\]
where the candidate library
\[
(\Theta_1,\ldots,\Theta_{12}) = \bigl(1,\ u,\ u^2,\ u^3,\ u_x,\ u_{xx},\ u_{xxx},\ u_{xxxx},\ u u_x,\ u^2 u_x,\ u u_{xx},\ u_x^2\bigr)
\]
contains polynomial, derivative, and interaction terms, and the target coefficient vector $\xi \in \mathbb{R}^M$ is sparse. The corresponding support, the set of active terms, is
\[
m_j = \mathbf{1}\{\xi_j \neq 0\} \in \{0,1\}.
\]

We call one space--time field $u(x,t)$ a \emph{trajectory} and a governing PDE with fixed coefficients a \emph{system}; trajectories of the same system differ only in their initial condition. A \emph{family} is an equation form with free coefficients (Appendix Table~\ref{tab:terminology}). A trajectory is observed with additive Gaussian noise,
\[
y = u + \sigma\,\mathrm{std}(u)\,\epsilon,
\qquad
\epsilon \sim \mathcal{N}(0,1),
\]
so a noise level of $\sigma = 0.8$, written $80\%$, means a noise standard deviation of $0.8$ times that of the field. Given one noisy trajectory $y$ together with its sampling metadata (grid spacing $\Delta x$, time step $\Delta t$, duration $T_{\mathrm{total}}$, domain length $L$, and the numbers of grid points and frames), the task is to identify both the support $m$ and the coefficients $\xi$.

\section{Method}
\label{sec:method}

UNIDE has three stages: denoise the observed trajectory (\S\ref{subsec:denoise}), select the active library terms (\S\ref{subsec:support}), and solve for coefficients on that support (\S\ref{subsec:coef}). The two learned stages share one denoiser and one encoder across all PDE families; the coefficient stage is a classical least-squares solve. Throughout, $y$ denotes the noisy observation of the clean field $u$ under the noise model of Section~\ref{sec:problem}.

\subsection{Self-Supervised Denoising}
\label{subsec:denoise}

Pointwise spatial derivatives amplify noise rapidly: the fourth-derivative channel can be corrupted by orders of magnitude even when the field noise is visually small. The first stage therefore removes the noise with a denoiser shared across PDE families,
\[
\tilde{u}=D_{\eta}(y),
\]
implemented as a residual dilated convolutional network with circular padding in space. It is trained with a blind-spot Noise2Void objective~\citep{krull2019noise2void}: a random mask $\Omega$ hides a small subset of space--time entries, $M_{\Omega}$ replaces each masked entry with the value of a randomly chosen nearby pixel so the network cannot copy the noise it must predict, and the network predicts the hidden values from their context. The loss is
\[
\mathcal{L}_{\mathrm{n2v}}(\eta)
=
\frac{1}{|\Omega|}
\sum_{(x,t)\in\Omega}
\left(D_{\eta}(M_{\Omega}y)_{x,t} - y_{x,t}\right)^2 .
\]
The loss compares the output only with the noisy observations, so the denoiser needs no clean targets, and at test time it requires no optimization, only one network evaluation per overlapping time window (Appendix~\ref{app:model}). At $10\%$ noise, removing the denoiser lowers Term F1 from $0.922$ to $0.861$ and raises coefficient error from $4.7\%$ to $14.6\%$, while exact match without pruning stays near $0.67$ (Table~\ref{tab:ablation-noise10}). Denoising thus mainly stabilizes the features used to select the terms; the coefficient fit below uses a weak form and avoids high-order pointwise differentiation of the field. Each family's space--time structure survives at $50\%$ noise (Appendix Figure~\ref{fig:denoise-fields}), and the residual relative $L_2$ error stays below the noise level everywhere and far below it from $5\%$ noise onward---from under $0.5\%$ at $1\%$ noise to about $5$--$8\%$ at $80\%$ (Appendix Figure~\ref{fig:denoise-summary}). Appendix~\ref{app:model} gives the architecture and training details.

\subsection{Term Selection on Denoised Features}
\label{subsec:support}

The second stage predicts which library terms are active. From the denoised field $\tilde{u}$ we compute the encoder input $X(\tilde{u})$, which stacks $u_t$ and the eleven non-constant library terms, together with time features $g_{1:T}$ and the sampling metadata $z$.

The encoder $E_\phi$ turns the sequence of $T$ frames into a single trajectory representation through a per-frame spatial model followed by a temporal model. Each frame $X_t$ is encoded by a spatial CNN,
\[
s_t = f_{\mathrm{cnn}}(X_t),
\]
each per-frame code is tagged with its time features,
\[
r_t = s_t + f_{\mathrm{time}}(g_t),
\]
and a transformer processes the tagged sequence across frames,
\[
h_{1:T} = f_{\mathrm{tr}}(r_{1:T}).
\]
Attention pooling over the frames, fused with an MLP embedding of the metadata $z$, gives the shared representation $h = E_\phi(X(\tilde{u}),g,z)$ (Appendix~\ref{app:model}). Linear heads read off $h$: a term head that predicts which terms are active, and auxiliary coefficient and family heads used only for the training losses.

The term head outputs support probabilities
\[
\ell = W_{\mathrm{term}}h + b_{\mathrm{term}},
\qquad
p=\sigma(\ell)\in[0,1]^M.
\]
At evaluation time, a validation-tuned threshold gives the hard support
\[
\hat{m}_j=\mathbf{1}\{p_j>\tau_{\mathrm{val}}\}.
\]
The threshold is tuned for validation F1 over the fixed grid $\{0.3, 0.4, \ldots, 0.9\}$; the selected value varies between $0.3$ and $0.8$ across noise levels, with no monotone trend. Training the encoder on denoised rather than clean features avoids a train--test mismatch that otherwise causes low-order terms to be over- or under-selected after denoising.

\subsection{Weak-Form Coefficient Identification}
\label{subsec:coef}

UNIDE does not predict coefficients with a neural network: once the support is predicted, they are fitted by weak-form least squares on that support, following WSINDy~\citep{messenger2021weak}. The PDE is multiplied by smooth, compactly supported test functions $\psi_k$ and integrated over local space--time subdomains. Integration by parts moves the time derivative and every pure spatial derivative, including fourth-order terms, from the data onto the test functions. Boundary terms vanish because the test functions and their required derivatives vanish on the subdomain boundaries. The implemented interaction columns for $u\,u_x$, $u^2u_x$, and $u\,u_{xx}$ retain first-order finite differences of the denoised field under the integrals (Appendix~\ref{app:weakform}). Thus, the weak coefficient assembly avoids high-order differentiation of the field, although the encoder still uses derivative features. This yields the linear system
\[
b_{\mathrm{w}} \approx \Theta_{\mathrm{w}}(\tilde{u})\,\xi,
\qquad
b_{\mathrm{w},k} = -\!\int (\partial_t \psi_k)\,\tilde{u},
\]
where each column of $\Theta_{\mathrm{w}}$ integrates a library term against the corresponding (differentiated) test functions. Given the predicted support $\hat{m}$, we solve
\[
\hat{\xi}
=
\arg\min_{\xi}
\left\|
\Theta_{\mathrm{w}}(\tilde{u})\operatorname{Diag}(\hat{m})\xi - b_{\mathrm{w}}
\right\|_2^2 .
\]
After the initial coefficient fit, greedy backward elimination removes a selected term when its deletion increases the normalized weak-form residual by no more than a validation-selected relative tolerance. The procedure reuses the assembled weak system but performs several least-squares fits per trajectory. It can remove excess terms but cannot introduce a missing term, and it can delete a true term with a small contribution to the fitted residual. Thus residual pruning is a model-selection heuristic on the observed trajectory, not a certificate of correct physics. Appendix~\ref{app:weakform} specifies the test functions, quadrature, normalized solves, and stopping rule.

\subsection{Training Objective}
\label{subsec:objective}

The encoder is trained with five loss terms:
\[
\mathcal{L}
=
\mathcal{L}_{\mathrm{term}}
+ \lambda_{\mathrm{coef}}\mathcal{L}_{\mathrm{coef}}
+ \lambda_{\mathrm{res}}\mathcal{L}_{\mathrm{res}}
+ \lambda_{\mathrm{sp}}\mathcal{L}_{\mathrm{sp}}
+ \lambda_{\mathrm{fam}}\mathcal{L}_{\mathrm{fam}}.
\]
The primary loss $\mathcal{L}_{\mathrm{term}}$ is a weighted binary cross-entropy between the predicted term probabilities and the true active terms; active terms are up-weighted because most library terms are inactive in any one equation. The remaining four losses are auxiliary: (1) the coefficient loss $\mathcal{L}_{\mathrm{coef}}$ is the relative absolute error between the auxiliary coefficient head and the true coefficients, evaluated on the active terms only; (2) the residual loss $\mathcal{L}_{\mathrm{res}}$ measures how well the library terms, combined with the auxiliary coefficients, reproduce the observed time derivative $u_t$, as the squared residual normalized by the magnitude of $u_t$; (3) the sparsity loss $\mathcal{L}_{\mathrm{sp}}$ is an $L_1$ penalty on the auxiliary coefficients, which favors equations with few terms; and (4) the family loss $\mathcal{L}_{\mathrm{fam}}$ is a cross-entropy loss for predicting the PDE family with the auxiliary family head.
The auxiliary coefficients serve only these training losses; every reported coefficient comes from the least-squares fit of Section~\ref{subsec:coef}. The denoiser is trained first, and the encoder is then trained on its output; both see per-sample noise drawn uniformly from $0$ up to the target level of the model (Section~\ref{sec:setup}). The loss weights ($\lambda_{\mathrm{coef}}{=}1$, $\lambda_{\mathrm{res}}{=}10^{-3}$, $\lambda_{\mathrm{sp}}{=}3{\times}10^{-2}$, $\lambda_{\mathrm{fam}}{=}10^{-1}$) are studied in Table~\ref{tab:loss-ablation}, and Appendix~\ref{app:model} gives the formal definitions and implementation details.

\section{Benchmark and Evaluation Setup}
\label{sec:dataset}

\begin{figure}[!htbp]
\centering
\includegraphics[width=\linewidth]{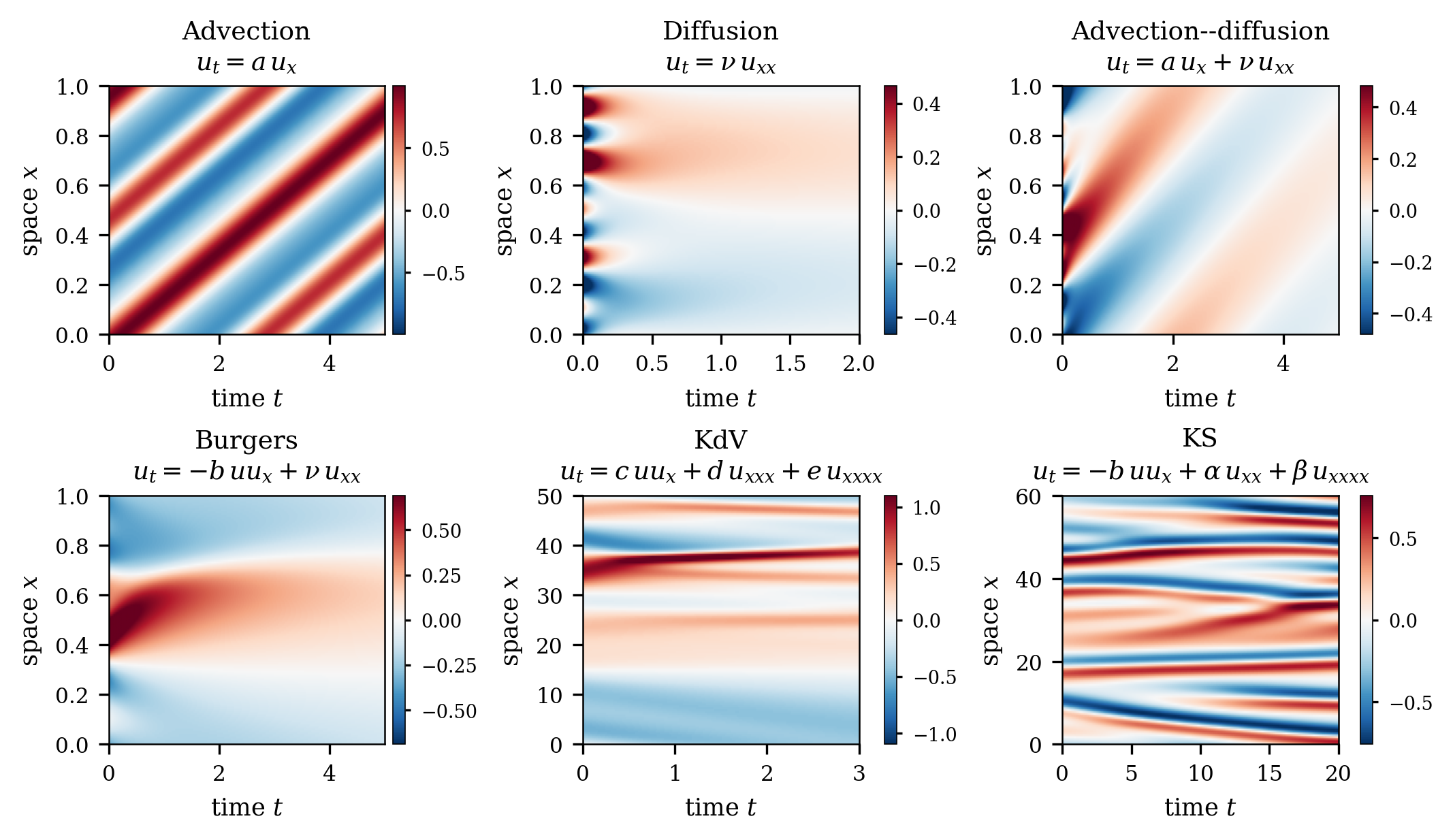}
\caption{A representative space--time solution $u(x,t)$ for each of the six PDE families, with its governing equation. The families span distinct dynamical regimes: advection translates the profile, diffusion smooths it, advection--diffusion does both, Burgers steepens then viscously regularizes, KdV produces dispersive waves, and Kuramoto--Sivashinsky is chaotic. Color shows $u$ (diverging scale, symmetric per panel); spatial domains differ by family.}
\label{fig:pde-gallery}
\end{figure}

\subsection{PDE Families and Data}

The benchmark contains six PDE families: advection, diffusion, advection--diffusion, Burgers, Korteweg--de Vries (KdV), and Kuramoto--Sivashinsky (KS). Each family contributes 20 coefficient settings with 20 trajectories each, giving 2,400 trajectories in total; every trajectory uses a distinct random initial condition, so no initial condition is shared across coefficient settings or families. The six families activate five of the 12 library terms overall: $u_x$, $u_{xx}$, $u_{xxx}$, $u_{xxxx}$, and $u u_x$. The benchmark tests identification across PDE families, coefficient scales, initial conditions, and observation-noise levels.

Trajectories are generated with pseudo-spectral solvers built on APEBench~\citep{apebench}, with equation-specific solver settings and subsampling. All stored trajectories use $N_x = 128$ spatial points and 501 time points on periodic domains; each sample includes the trajectory, sampling metadata, and precomputed library channels. Appendix~\ref{app:data} lists the solver settings, sampling grids, and initial-condition distribution. The families differ in both their active terms and their coefficient ranges:
\begin{itemize}
    \item \textbf{Advection:} $u_t = a\,u_x$. This is the pure transport family, with only the $u_x$ term active. The advection coefficient is sampled from $(-1.0,-0.05)\cup(0.05,1.0)$ so that both directions of transport are represented while excluding near-zero drift.
    \item \textbf{Diffusion:} $u_t = \nu\,u_{xx}$. Only the $u_{xx}$ term is active. The diffusion coefficient is sampled from $(0.003,0.015)$, calibrated so trajectories visibly smooth over the observed window without collapsing immediately to zero.
    \item \textbf{Advection--diffusion:} $u_t = a\,u_x + \nu\,u_{xx}$. The active terms are $u_x$ and $u_{xx}$. Coefficients are sampled as $a\in(-0.5,-0.1)$ and $\nu\in(0.003,0.012)$ so that both transport and smoothing remain visible.
    \item \textbf{Burgers:} $u_t = -b\,u\,u_x + \nu\,u_{xx}$. The active terms are $u\,u_x$ and $u_{xx}$. Coefficients are sampled as $b\in(0.1,0.5)$ and $\nu\in(0.001,0.005)$, producing an interplay between nonlinear steepening and viscous smoothing.
    \item \textbf{KdV:} $u_t = c\,u\,u_x + d\,u_{xxx} + e\,u_{xxxx}$. We use a stabilized KdV variant with a fixed fourth-derivative term, which is included in the ground-truth support and coefficient evaluation. The active terms are $u\,u_x$, $u_{xxx}$, and $u_{xxxx}$. Coefficients vary over $c\in(-8.0,-2.0)$ and $d\in(-2.0,-0.5)$, while $e=-0.125$ is fixed for numerical stabilization.
    \item \textbf{KS:} $u_t = -b\,u\,u_x + \alpha\,u_{xx} + \beta\,u_{xxxx}$. The active terms are $u\,u_x$, $u_{xx}$, and $u_{xxxx}$. Coefficients are sampled over $b\in(1.0,4.0)$, $\alpha\in(-1.2,-0.5)$, and $\beta\in(-0.8,-0.4)$ so that long-wave instability and short-wave damping are both present.
\end{itemize}
The time window of each family is chosen so that the stored 500-step trajectory resolves the dynamics accurately while still showing visible temporal evolution (Figure~\ref{fig:pde-gallery}).

\subsection{Test Sets and Observation Noise}
\label{sec:testsets}

We evaluate UNIDE on two test sets. The random-trajectory test set contains new noisy trajectories of the systems seen in training; the coefficient-extrapolation test set contains trajectories of unseen systems with larger coefficient-vector norms. Clean versions of both serve as references.

\paragraph{Random-trajectory generalization.}
Within each PDE family, trajectories are partitioned at random into 200 / 50 / 150 train / validation / test samples, so the full test set contains 900 trajectories. Because the partition is over trajectories, each test trajectory is a new initial condition of one of the 120 systems seen in training; its coefficients are still estimated by least squares from the trajectory itself, and unseen systems are tested by the coefficient-extrapolation setting below. In the noisy version---the paper's main benchmark---Gaussian noise is added to the test trajectories and every method must identify the governing equation from the corrupted observations; the clean version is reported as a reference in Appendix~\ref{app:clean}.

\paragraph{Coefficient extrapolation.}
Within each PDE family, the 20 coefficient settings are sorted by the norm of their active coefficients; the 14 smallest are used for training, the next 3 for validation, and the 3 largest for testing, so every test system has a larger coefficient norm than any training system. We report this setting at 10\% observation noise, where distribution shift and derivative corruption act together; clean extrapolation results appear in Appendix~\ref{app:clean}. A separate denoiser and encoder are trained only on the training partition of this setting, so no test system is seen during training.

\paragraph{Observation noise.}
The noise level takes the values $\sigma \in \{0,\allowbreak 0.01,\allowbreak 0.05,\allowbreak 0.10,\allowbreak 0.20,\allowbreak \ldots,\allowbreak 0.80\}$ in the observation model of Section~\ref{sec:problem}. All derivative features and weak-form systems are recomputed from the noisy observations. The random-trajectory noise sweep is the main aggregate result in Section~\ref{sec:results}, while the 10\% coefficient-extrapolation result serves as a targeted stress test.

\subsection{Baselines and Training Protocol}
\label{sec:setup}

We compare against three baselines that span the main discovery paradigms. Each is evaluated per trajectory on the same stored fields, with sparse-regression thresholds tuned on noisy validation trajectories at each noise level:
\begin{itemize}
    \item \textbf{SINDy}~\citep{brunton2016discovering,rudy2017data}: pointwise (strong-form) sparse regression with sequentially thresholded least squares (STLSQ) on finite-difference derivative features. In the PDE setting, it is equivalent to PDE-FIND.
    \item \textbf{WSINDy}~\citep{messenger2021weak}: weak-form sparse regression, with coefficients re-solved by weak-form least squares on its own selected support so that the baseline is not penalized by a noise-fragile pointwise refit.
    \item \textbf{PySR}~\citep{cranmer2023pysr}: per-trajectory genetic-programming symbolic regression over the library features, followed by a least-squares refit of the coefficients of the library terms that appear in its expression. The refit is needed because PySR's expression can combine library terms nonlinearly and its constants are fitted on a 4,000-point subsample.
\end{itemize}
Appendix~\ref{app:eval} reports their candidate libraries, output mappings, validation subsets, and hyperparameters. Because these pipelines differ from UNIDE in preprocessing and tuning budget, Section~\ref{sec:matched-controls} adds controls that give sparse regression the same denoiser, candidate terms, pruning, and validation data as UNIDE.

The main noise sweep uses a separate UNIDE model for each target observation-noise level $\sigma$, which we call a specialist. Each specialist comprises a denoiser and an encoder trained sequentially across the six PDE families, using per-sample noise drawn from $U(0,\sigma)$. The baselines fit each noisy trajectory directly and have no training range; the PySR sweep covers the full 900-trajectory test set at every reported level. Appendix~\ref{app:perlevel} evaluates transfer across noise levels.

\subsection{Metrics}

An identified equation can fail in two ways: it can contain the wrong terms, or the right terms with the wrong coefficients. We therefore report two structure metrics and two coefficient metrics. \textbf{Exact match}, our primary metric, measures the fraction of test trajectories whose predicted active term set equals the true one, irrespective of coefficient accuracy. \textbf{Term F1} asks how close the term set is when recovery is not exact, averaging the per-trajectory F1 of the selected terms. \textbf{Coefficient RMSE} and \textbf{coefficient relative error} ask how accurately the coefficients of the true terms are estimated; both are computed only on the true active terms. The two families are complementary by construction: a spurious extra term lowers structure but not coefficient error, while a missing true term lowers both, so the structure and coefficient columns of every table should be read together. Formal definitions and a worked example are given in Appendix~\ref{app:eval}.

\section{Results}
\label{sec:results}

\subsection{Random-Trajectory Generalization under Observation Noise}

\begin{figure}[!htbp]
\centering
\includegraphics[width=\linewidth]{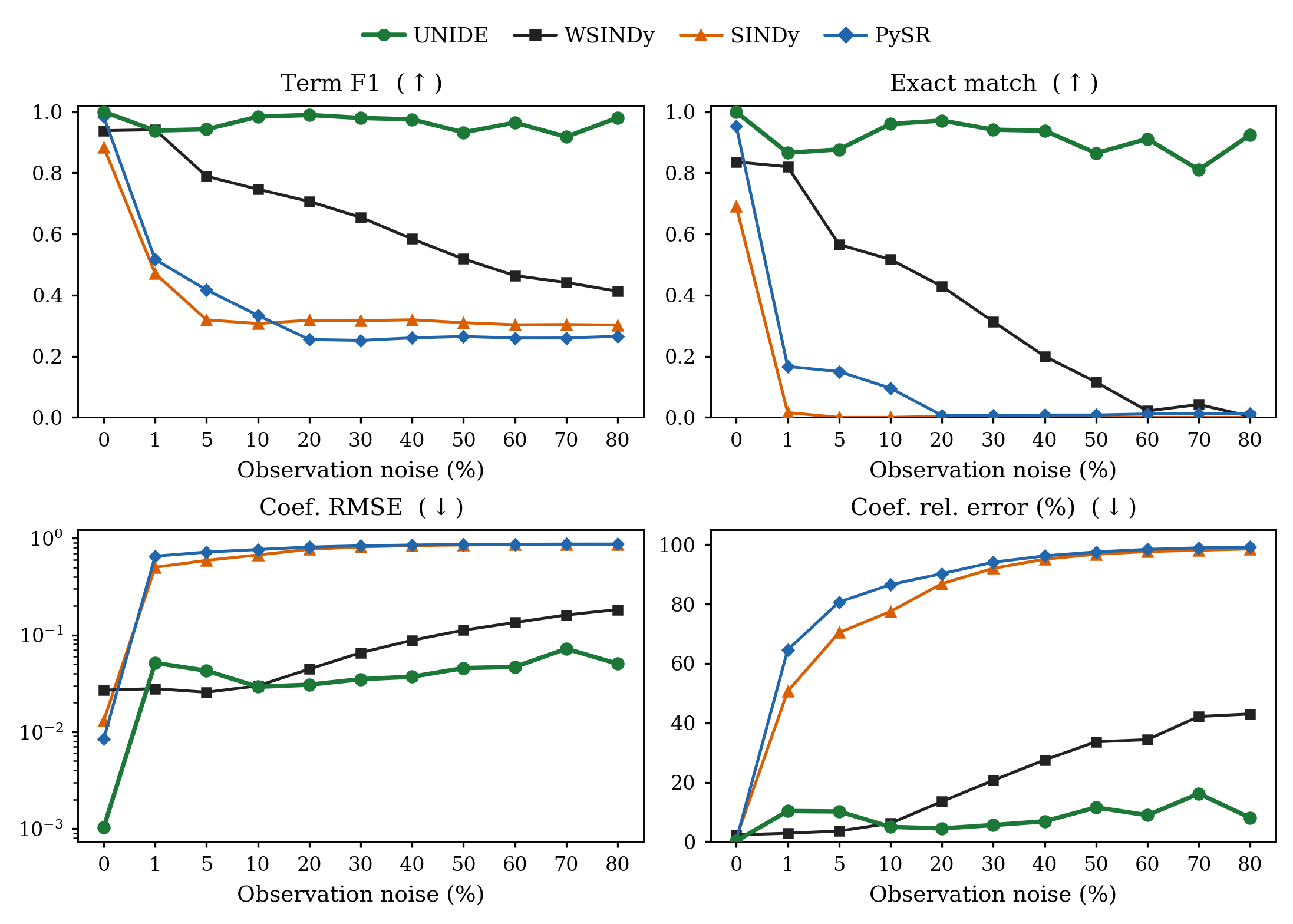}
\caption{Noise sweep on the random-trajectory test set (900 held-out trajectories from six PDE families). The four panels report exact match, Term F1, coefficient RMSE, and coefficient relative error against the observation-noise level. UNIDE uses the specialist trained at each level and includes its denoiser; SINDy, WSINDy, and PySR receive the raw noisy fields, with the candidate libraries detailed in Appendix~\ref{app:baselines}. UNIDE has the highest exact match at every noisy level. At $0\%$ there is no noise to remove, so UNIDE skips the denoiser and fits coefficients by least squares on the stored field (Appendix Table~\ref{tab:clean-ref}).}
\label{fig:main-results}
\end{figure}

\begin{figure}[!htbp]
\centering
\includegraphics[width=\linewidth]{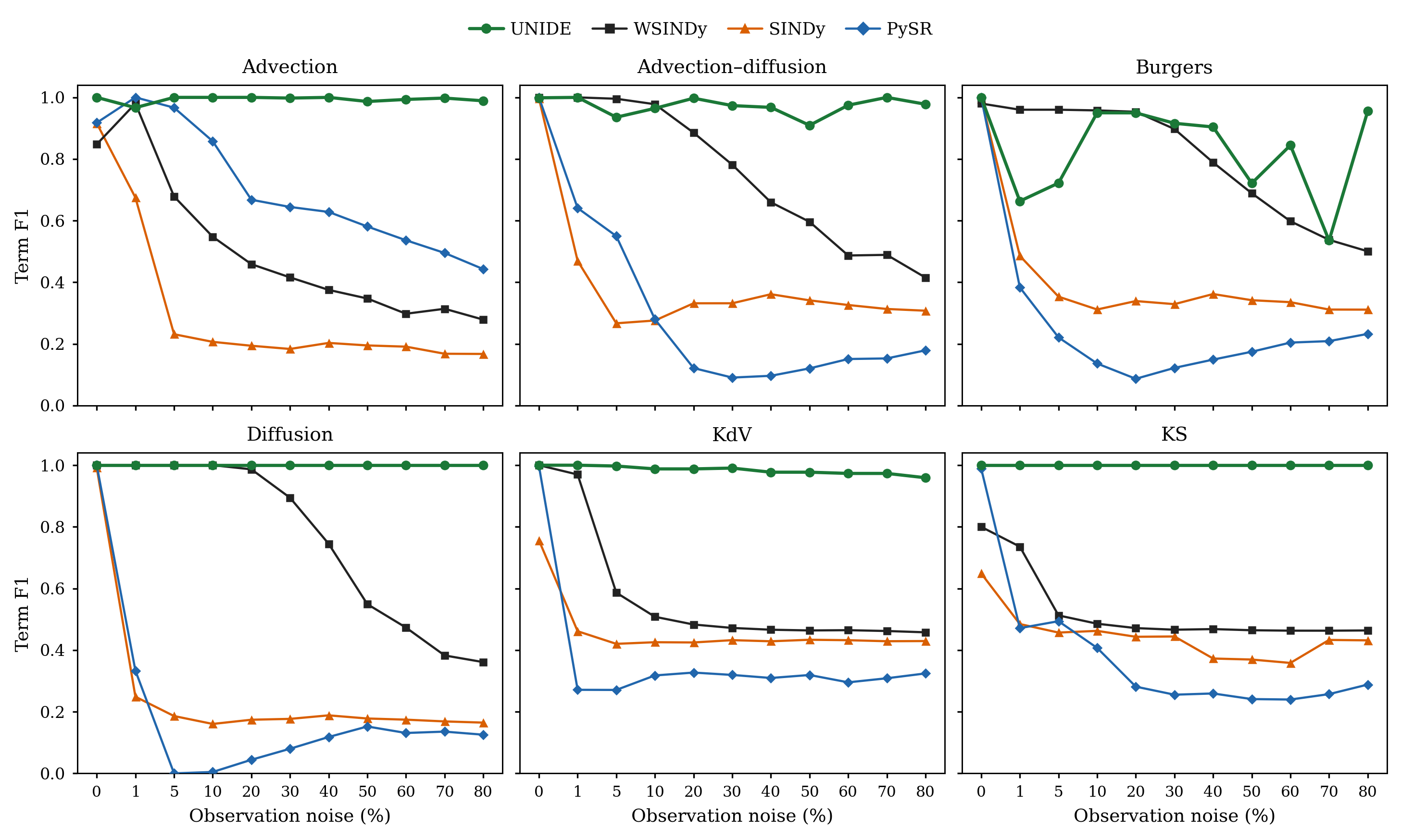}
\caption{Term F1 versus observation-noise level for each PDE family on the random-trajectory test set (150 trajectories per family; UNIDE trained at each noise level, with the pruned solve). UNIDE keeps Term F1 above $0.90$ on five of the six families across the whole sweep; Burgers is the exception. At $0\%$, UNIDE skips the denoiser (Appendix Table~\ref{tab:clean-ref}). The per-family exact match, coefficient RMSE, and relative error are in Appendix Figures~\ref{fig:per-pde-exact}--\ref{fig:per-pde-relerr}.}
\label{fig:per-pde}
\end{figure}

In this section, we evaluate the central test of this study: identifying the governing equation from noisy observations. We sweep the observation-noise level and evaluate all methods on the same 900 held-out trajectories. Figure~\ref{fig:main-results} summarizes the aggregate trend, Appendix Table~\ref{tab:noise-sweep} lists its values, and Figure~\ref{fig:per-pde} breaks it down by family. The reported performance combines three design choices: training on 1,200 trajectories from the six families, matching each specialist's training-noise range to its evaluation level, and applying residual pruning after term selection.

\paragraph{Term recovery.}
UNIDE has the highest exact match at every noisy level, and the margin widens with noise. At 80\% noise UNIDE's exact match is 0.924, against 0.003 for WSINDy, and the pointwise methods collapse outright. Term F1 shows the same trend: UNIDE has the highest Term F1 from $5\%$ onward, and between 5\% and 80\% noise its Term F1 stays above 0.91 (0.942 at 5\%, 0.980 at 80\%), while WSINDy falls from 0.789 to 0.413; at 1\% the two have similar Term F1 (0.938 versus 0.941). The gap opens only under noise: on clean data, the per-trajectory baselines are already strong, and PySR attains the highest clean exact match of the baselines and, on the random-trajectory reference, their lowest relative coefficient error (Appendix~\ref{app:clean}).

\paragraph{Coefficient accuracy.}
The advantage in coefficient accuracy appears only at higher noise: up to 5\% noise WSINDy's weak-form solve is more accurate than UNIDE's; at 10\% UNIDE's relative error is lower (5.0\% versus 6.2\%), but the five-seed mean is not (7.3\%); and from 20\% onward UNIDE attains the lowest coefficient error by a wide margin (8.0\% versus 43.0\% relative error at 80\%).

\paragraph{Seed variability.}
The reported encoder seed is seed 0, the training default for every run, fixed before the seed study. Across five encoder seeds, the standard deviation of exact match is at most $0.052$, and seed 0 lies within $1.6$ standard deviations of the seed mean on every metric; it is the strongest of the five at 10\% and 80\% and the weakest at 50\% (Appendix~\ref{app:seeds}).

\paragraph{Per-family results.}
The aggregate hides per-family structure. Because the test set is balanced (150 trajectories per family), the aggregate metrics are an unweighted macro-average over the six families. With the pruned solve, advection, diffusion, and KS hold exact match between 0.96 and 1.0 across the entire sweep, and advection--diffusion between 0.71 and 1.0. SINDy fails early because pointwise derivative features are corrupted; WSINDy is more stable on low-order dissipative equations, but its Term F1 stalls near $0.5$ on the dispersive and chaotic families from $10\%$ noise onward.

Two failure modes remain. For Burgers, UNIDE's exact match is 0.23 at 1\% noise, 0.47 at 5\%, and 0.00 at 70\%, against 0.63--0.93 at the other levels (Appendix Figure~\ref{fig:per-pde-exact}). Burgers differs from advection--diffusion only in its transport term, $u\,u_x$ instead of $u_x$, so a single wrong term turns one family into the other. Most Burgers errors are exactly this swap, which pruning cannot repair because it only removes terms. The second failure is on KdV, where pruning removes a true term. KdV contains a small fixed term, $-0.125\,u_{xxxx}$, added for numerical stability, and deleting it changes the weak-form fit so little that the change can fall within the pruning tolerance. Pruning then removes it on up to 29 of 150 trajectories, which keeps KdV exact match between 0.81 and 1.0.

\FloatBarrier

\subsection{Learned Term Selection under Matched Denoising}
\label{sec:matched-controls}

The comparison above gives UNIDE a denoiser that the baselines lack, so its advantage could come from denoising alone. We test this with controls that give weak-form sparse regression (STLSQ) everything UNIDE uses except learned term selection: the same denoised field, the same ten candidate terms (the weak-form library omits $1$ and $u_x^2$, which no family contains) and 400 weak-form windows, the same residual pruning, and hyperparameters selected on the same full validation partition (Tables~\ref{tab:audit-controls-A10}--\ref{tab:audit-controls-B10}; settings in Appendix~\ref{app:baselines}). Learned term selection has the higher exact match in every setting: $0.960$ compared with $0.787$ at $10\%$ noise, $0.924$ versus $0.834$ at $80\%$ noise (five-seed mean $0.888$), and $0.906$ compared with $0.642$ for held-out systems with larger coefficients at $10\%$ noise (Section~\ref{sec:extrap}). Denoising and pruning each help sparse regression but do not close the gap: at $80\%$ noise its exact match rises from $0.010$ on the raw field to $0.484$ after denoising and $0.834$ with pruning.

Coefficient accuracy shows no consistent winner. Fitting the same weak-form system on its own selected support, the matched sparse selector has lower mean relative coefficient error at $10\%$ noise ($3.28\%$ versus $4.99\%$) and under coefficient extrapolation ($3.27\%$ versus $9.88\%$), and higher error at $80\%$ noise ($10.14\%$ versus $8.04\%$). Learned selection thus improves which terms are found rather than how accurately their coefficients are fitted. Assigning each trajectory the fixed support of the family predicted by the encoder's auxiliary head reaches exact match $0.760$--$0.847$ across the three settings, with or without pruning; this diagnostic uses the auxiliary head, not a separately optimized family classifier.

% Generated from completed experiment outputs only.
\begin{table}[!htbp]
\centering\small
\setlength{\tabcolsep}{3pt}
\caption{Matched controls: random-trajectory test set, 10\% noise (900 test trajectories). Every weak-form row uses the same ten candidate terms and $K=400$ windows, and the denoised rows share UNIDE's denoiser. Hyperparameters are selected on the full validation partition. The auxiliary family-head rows assign each trajectory the fixed term set of the family predicted by the encoder's auxiliary head. Coefficient error is evaluated on the true active terms. Best per metric in bold.}
\label{tab:audit-controls-A10}
\begin{tabularx}{\linewidth}{@{}Xrrrr@{}}\toprule
Method & Exact $\uparrow$ & Term F1 $\uparrow$ & Coef. RMSE $\downarrow$ & Coef. rel\% $\downarrow$ \\ \midrule
UNIDE & \textbf{0.960} & \textbf{0.984} & 2.92e-02 & 4.99 \\
Weak STLSQ, raw field & 0.484 & 0.753 & 3.42e-02 & 8.23 \\
Denoiser + weak STLSQ & 0.504 & 0.770 & 3.41e-02 & 9.13 \\
Denoiser + weak STLSQ + pruning & 0.787 & 0.946 & \textbf{2.53e-02} & \textbf{3.28} \\
Denoiser + all-term backward pruning & 0.402 & 0.833 & 2.54e-02 & 3.99 \\
Auxiliary family head & 0.760 & 0.866 & 5.67e-02 & 18.53 \\
Auxiliary family head + pruning & 0.761 & 0.869 & 5.67e-02 & 18.56 \\
\bottomrule\end{tabularx}
\end{table}

\begin{table}[!htbp]
\centering\small
\setlength{\tabcolsep}{3pt}
\caption{Matched controls: random-trajectory test set, 80\% noise (900 test trajectories). Rows and settings as in Table~\ref{tab:audit-controls-A10}.}
\label{tab:audit-controls-A80}
\begin{tabularx}{\linewidth}{@{}Xrrrr@{}}\toprule
Method & Exact $\uparrow$ & Term F1 $\uparrow$ & Coef. RMSE $\downarrow$ & Coef. rel\% $\downarrow$ \\ \midrule
UNIDE & \textbf{0.924} & \textbf{0.980} & \textbf{5.04e-02} & \textbf{8.04} \\
Weak STLSQ, raw field & 0.010 & 0.439 & 1.82e-01 & 42.64 \\
Denoiser + weak STLSQ & 0.484 & 0.712 & 6.28e-02 & 10.06 \\
Denoiser + weak STLSQ + pruning & 0.834 & 0.954 & 5.90e-02 & 10.14 \\
Denoiser + all-term backward pruning & 0.613 & 0.897 & 5.74e-02 & 9.61 \\
Auxiliary family head & 0.760 & 0.856 & 1.05e-01 & 21.24 \\
Auxiliary family head + pruning & 0.767 & 0.863 & 1.05e-01 & 21.44 \\
\bottomrule\end{tabularx}
\end{table}

\begin{table}[!htbp]
\centering\small
\setlength{\tabcolsep}{3pt}
\caption{Matched controls: coefficient-extrapolation test set, 10\% noise (360 test trajectories). Rows and settings as in Table~\ref{tab:audit-controls-A10}.}
\label{tab:audit-controls-B10}
\begin{tabularx}{\linewidth}{@{}Xrrrr@{}}\toprule
Method & Exact $\uparrow$ & Term F1 $\uparrow$ & Coef. RMSE $\downarrow$ & Coef. rel\% $\downarrow$ \\ \midrule
UNIDE & \textbf{0.906} & \textbf{0.948} & 8.79e-02 & 9.88 \\
Weak STLSQ, raw field & 0.472 & 0.738 & \textbf{4.44e-02} & 5.98 \\
Denoiser + weak STLSQ & 0.475 & 0.734 & 4.60e-02 & 5.93 \\
Denoiser + weak STLSQ + pruning & 0.642 & 0.898 & 4.58e-02 & \textbf{3.27} \\
Denoiser + all-term backward pruning & 0.317 & 0.789 & 4.55e-02 & 4.14 \\
Auxiliary family head & 0.833 & 0.877 & 1.43e-01 & 15.83 \\
Auxiliary family head + pruning & 0.847 & 0.887 & 1.43e-01 & 15.83 \\
\bottomrule\end{tabularx}
\end{table}

\FloatBarrier

\subsection{Ablation Study at 10\% Observation Noise}

We ablate the pipeline at 10\% observation noise on the full 900-trajectory test set. Table~\ref{tab:ablation-noise10} isolates the four coupled design choices (denoising the observed field, training the encoder on denoised features, estimating coefficients with a weak-form solve, and pruning the predicted support using the weak-form residual) and, in its last row, replaces the term head with the encoder's auxiliary family head and fixed support templates as a diagnostic comparison. Table~\ref{tab:loss-ablation} then removes one auxiliary loss term at a time.

\begin{table}[!htbp]
\centering
\small
\setlength{\tabcolsep}{4pt}
\caption{Ablations of the pipeline at 10\% observation noise on the full 900-trajectory test set. Each row changes one component of UNIDE; the last row replaces the term head with the encoder's family-classification head and fixed per-family term sets. The Full UNIDE row is the $10\%$ point of the noise sweep (Figure~\ref{fig:main-results}). Rows other than the first are evaluated without residual pruning, so the ``No residual pruning'' row is their shared reference. Best per metric in bold.}
\label{tab:ablation-noise10}
\begin{tabular}{lcccc}
\toprule
Variant & Exact $\uparrow$ & Term F1 $\uparrow$ & Coef. RMSE $\downarrow$ & Coef. rel\% $\downarrow$ \\
\midrule
\textbf{Full UNIDE} & \textbf{0.960} & \textbf{0.984} & 2.92e-02 & 4.99 \\
No residual pruning & 0.670 & 0.922 & \textbf{2.87e-02} & \textbf{4.73} \\
No denoiser & 0.667 & 0.861 & 5.74e-02 & 14.60 \\
Clean-trained encoder & 0.393 & 0.603 & 1.51e-01 & 38.51 \\
Strong-form coefficient LS & 0.670 & 0.922 & 5.11e-01 & 48.85 \\
Family-template support & 0.760 & 0.866 & 5.67e-02 & 18.53 \\
\bottomrule
\end{tabular}
\end{table}

At $10\%$ noise, removing residual pruning lowers exact match from $0.960$ to $0.670$. Relative coefficient error changes from $4.99\%$ to $4.73\%$, so pruning trades a small increase in coefficient error for the gain in exact match. Against the unpruned reference, using the clean-trained encoder with the denoised fields gives F1 $0.603$ and relative coefficient error $38.5\%$, while removing denoising gives F1 $0.861$ and error $14.6\%$. Replacing weak-form fitting with pointwise least squares gives error $48.9\%$. The matched controls (Section~\ref{sec:matched-controls}) show that a simpler denoising-plus-sparse-regression pipeline does not reach the same term recovery.

The last row assigns each trajectory the fixed support of the family predicted by the encoder's auxiliary head. Without pruning, it attains exact match $0.760$ versus $0.670$ for the term head, but F1 $0.866$ versus $0.922$ and relative coefficient error $18.5\%$ versus $4.7\%$, because each of its errors replaces the whole support with another family's. With pruning, the term head reaches exact match $0.960$, whereas the family head improves only to $0.761$ (Table~\ref{tab:audit-controls-A10}).

\paragraph{Per-loss-term ablation.}
Table~\ref{tab:loss-ablation} removes one auxiliary loss at a time, retraining the encoder from the same warm start with the same data, seed, and denoiser. The primary term loss is always kept, and every row uses the unpruned solve. Exact match ranges from $0.679$ to $0.711$ and relative coefficient error from $4.6\%$ to $6.3\%$, compared with $0.670$ and $4.7\%$ for the full objective. Every removal leaves exact match at or above that of the full objective, so in these single runs no auxiliary loss is needed for term recovery; removing the family or residual loss does raise coefficient error ($6.3\%$ and $5.3\%$ versus $4.7\%$).

\begin{table}[!htbp]
\centering
\small
\setlength{\tabcolsep}{4pt}
\caption{Per-loss-term ablation at 10\% observation noise on all 900 trajectories of the random-trajectory test set. Each row retrains the encoder with one auxiliary loss term removed, holding the denoiser, data, seed, and all other loss weights fixed; the primary term loss $\mathcal{L}_{\mathrm{term}}$ is always kept. Every row uses the unpruned solve and is a single run. Best per metric in bold.}
\label{tab:loss-ablation}
\begin{tabular}{lcccc}
\toprule
Objective & Exact $\uparrow$ & Term F1 $\uparrow$ & Coef. RMSE $\downarrow$ & Coef. rel\% $\downarrow$ \\
\midrule
\textbf{Full objective} & 0.670 & 0.922 & \textbf{2.87e-02} & 4.73 \\
$-\,\lambda_{\mathrm{coef}}$ (coefficient MAE) & 0.680 & \textbf{0.925} & 3.03e-02 & 4.83 \\
$-\,\lambda_{\mathrm{res}}$ (residual consistency) & 0.698 & 0.923 & 3.30e-02 & 5.30 \\
$-\,\lambda_{\mathrm{sp}}$ (L1 sparsity) & 0.679 & 0.924 & 2.99e-02 & \textbf{4.59} \\
$-\,\lambda_{\mathrm{fam}}$ (family CE) & \textbf{0.711} & 0.922 & 3.49e-02 & 6.34 \\
\bottomrule
\end{tabular}
\end{table}

\FloatBarrier

\subsection{Coefficient Extrapolation under 10\% Observation Noise}
\label{sec:extrap}

Coefficient extrapolation asks whether identification extends to systems never seen in training (Section~\ref{sec:testsets}). The coefficient settings are unseen, while the active term sets are represented in training; individual coefficients can still fall within their training ranges. Table~\ref{tab:coef-extrap-noise10} compares the methods on these held-out systems.

\begin{table}[!htbp]
\centering
\small
\setlength{\tabcolsep}{4pt}
\caption{Coefficient extrapolation under 10\% observation noise. Test trajectories come from held-out large-coefficient settings; UNIDE uses a denoiser and encoder trained only on the coefficient-extrapolation training partition.}
\label{tab:coef-extrap-noise10}
\begin{tabular}{lcccc}
\toprule
Method & Exact $\uparrow$ & Term F1 $\uparrow$ & Coef. RMSE $\downarrow$ & Coef. rel\% $\downarrow$ \\
\midrule
SINDy & 0.000 & 0.310 & 9.76e-01 & 77.87 \\
WSINDy & 0.494 & 0.717 & \textbf{4.37e-02} & \textbf{4.98} \\
PySR & 0.019 & 0.377 & 1.10e+00 & 81.23 \\
\textbf{UNIDE} & \textbf{0.906} & \textbf{0.948} & 8.79e-02 & 9.88 \\
\bottomrule
\end{tabular}
\end{table}

Under simultaneous extrapolation and 10\% noise, UNIDE recovers structure more accurately than the baselines (exact match 0.906 and Term F1 0.948, versus 0.494 and 0.717 for WSINDy), but its coefficients are less accurate than WSINDy's (relative error 9.9\% versus 5.0\%), largely because of Burgers, where the encoder again substitutes $u_x$ for $u\,u_x$ (exact match 0.53 and 33\% relative error, while the other five families reach exact match of at least 0.91; Appendix Table~\ref{tab:app_per_pde_noise10}). The pointwise methods fail as before. The advantage in term recovery therefore extends to held-out systems with larger coefficient-vector norms, while the coefficients remain less accurate than WSINDy's.

\FloatBarrier

\subsection{Inference Cost}
\label{sec:cost}

\begin{table}[!htbp]
\centering
\small
\setlength{\tabcolsep}{4pt}
\caption{Inference cost per trajectory at 10\% observation noise (hardware in Appendix~\ref{app:model}). All methods are timed on the same 18 trajectories on the same idle machine, with the full pipeline of each method including UNIDE's residual pruning; each entry is the mean $\pm$ standard deviation over five repeated runs after a discarded warm-up (threshold tuning is not counted, as it is not a per-trajectory cost). Exact match is measured at the same $10\%$ noise level (Appendix Table~\ref{tab:noise-sweep}); SINDy is faster but fails under noise.}
\label{tab:efficiency}
\begin{tabular}{lcccl}
\toprule
Method & ms / traj & Relative to UNIDE & Exact match & What runs \\
\midrule
SINDy & $43.4 \pm 0.4$ & $0.1\times$ & 0.000 & Pointwise sparse regression + STLSQ \\
\textbf{UNIDE} & $\mathbf{451.5 \pm 1.4}$ & $\mathbf{1\times}$ & \textbf{0.960} & Denoiser + encoder + pruned weak-form LS \\
WSINDy & $1198.0 \pm 3.4$ & $2.7\times$ & 0.517 & Weak-form search + coefficients \\
PySR & $22924.4 \pm 82.1$ & $51\times$ & 0.096 & Genetic symbolic search \\
\bottomrule
\end{tabular}
\end{table}

In the measured timing comparison, UNIDE is faster than both search-based methods: about $2.7\times$ faster than WSINDy and $51\times$ faster than PySR per trajectory; only pointwise SINDy is faster, and it fails under noise. The learned stages take only $56$~ms (denoiser $48$~ms, encoder $8$~ms), and the remaining $395$~ms is the weak-form coefficient solve, including the residual pruning.

\FloatBarrier

\section{Discussion}

The proposed UNIDE separates the two halves of equation identification: learned components decide which terms are present, and a least-squares fit decides their coefficients. Unlike per-trajectory sparse and symbolic regression, which fit each noisy field from scratch, UNIDE reuses what it learns across PDE families, so a new trajectory needs only windowed denoising and one encoder evaluation, followed by a few small least-squares solves.

We provide three reasons for its performance under heavy noise. (1) Learned term selection: given the same denoised field, candidate terms, and pruning, learned term selection recovers more exact term sets than sparse regression in every tested setting, whereas neither method fits coefficients consistently better (Section~\ref{sec:matched-controls}). (2) Training on denoised features: the encoder is trained on the output of the denoiser it sees at test time; a clean-trained encoder applied to the same denoised fields lowers exact match from $0.670$ to $0.393$ at $10\%$ noise (Table~\ref{tab:ablation-noise10}). (3) Weak-form coefficient fitting: the weak form moves derivatives from the data onto smooth test functions; replacing it with pointwise least squares raises the coefficient error from $4.7\%$ to $48.9\%$ (Table~\ref{tab:ablation-noise10}).

Two systematic failures remain beneath the aggregate: on Burgers the encoder can substitute $u_x$ for $u\,u_x$, which removal-only pruning cannot repair, and pruning can delete KdV's small stabilizing term. The first could be addressed by letting the solve swap terms as well as remove them.

\paragraph{Limitations.}
The benchmark is built to test identification under observation noise and coefficient shift, and its design sets the scope of these conclusions. All fields are synthetic, one-dimensional, periodic, and fully observed, and every test equation shares one of the six term structures seen in training. With known structures, choosing among family templates is a strong alternative to term-level prediction, and our experiments do not show that term-level prediction is necessary in this setting; held-out operator combinations are the test that separates the two. Each noise level uses its own specialist; with weights and thresholds frozen at their $10\%$ settings, exact match falls to $0.489$ at $80\%$ noise (Appendix Table~\ref{tab:fully-frozen}). Seed variability is measured for the encoder with the denoiser and warm start fixed, not for the full pipeline.

Future work may include (1) held-out operator combinations and held-out PDE families, which test whether the encoder has learned individual operators rather than family templates; (2) physical measurements; and (3) higher-dimensional domains with their candidate libraries, which require broader training data. Together these experiments lead from the present results toward a foundation model that discovers governing equations in new applications.

\section{Conclusion}

In this study, we proposed UNIDE, a data-driven approach that learns across PDE families to decide which terms govern a noisy trajectory and fits their coefficients by weak-form least squares. On six one-dimensional PDE families, it recovers the exact term set more often than SINDy, WSINDy, and PySR at every noise level from $1\%$ to $80\%$. At $80\%$ noise it recovers the exact term set for $92\%$ of 900 held-out trajectories (five-seed mean $89\%$), and it identifies each trajectory $51$ times faster than PySR. The advantage comes from learned term selection: given the same denoised field, candidate terms, and pruning, sparse regression recovers fewer exact term sets in every setting, including held-out systems with larger coefficients ($64\%$ versus $91\%$ for UNIDE at $10\%$ noise). Neither method fits the coefficients consistently more accurately. Extending these results from known term structures to new operator combinations and physical measurements is the next step toward a model that discovers equations in new applications.

\bibliographystyle{unsrtnat}
\bibliography{references}

\begin{thebibliography}{18}
\providecommand{\natexlab}[1]{#1}
\providecommand{\url}[1]{\texttt{#1}}
\expandafter\ifx\csname urlstyle\endcsname\relax
  \providecommand{\doi}[1]{doi: #1}\else
  \providecommand{\doi}{doi: \begingroup \urlstyle{rm}\Url}\fi

\bibitem[Brunton et~al.(2016)Brunton, Proctor, and
  Kutz]{brunton2016discovering}
Steven~L. Brunton, Joshua~L. Proctor, and J.~Nathan Kutz.
\newblock Discovering governing equations from data by sparse identification of
  nonlinear dynamical systems.
\newblock \emph{Proceedings of the National Academy of Sciences}, 113\penalty0
  (15):\penalty0 3932--3937, 2016.

\bibitem[Rudy et~al.(2017)Rudy, Brunton, Proctor, and Kutz]{rudy2017data}
Samuel~H. Rudy, Steven~L. Brunton, Joshua~L. Proctor, and J.~Nathan Kutz.
\newblock Data-driven discovery of partial differential equations.
\newblock \emph{Science Advances}, 3\penalty0 (4):\penalty0 e1602614, 2017.

\bibitem[Xie et~al.(2022)Xie, Samaei, Guo, Liu, and Gan]{xie2022data}
Xiaoyu Xie, Arash Samaei, Jiachen Guo, Wing~Kam Liu, and Zhengtao Gan.
\newblock Data-driven discovery of dimensionless numbers and governing laws
  from scarce measurements.
\newblock \emph{Nature Communications}, 13\penalty0 (1):\penalty0 7562, 2022.

\bibitem[Cranmer(2023)]{cranmer2023pysr}
Miles Cranmer.
\newblock Interpretable machine learning for science with {PySR} and
  {SymbolicRegression.jl}.
\newblock \emph{arXiv preprint arXiv:2305.01582}, 2023.

\bibitem[Messenger and Bortz(2021)]{messenger2021weak}
Daniel~A. Messenger and David~M. Bortz.
\newblock Weak {SINDy} for partial differential equations.
\newblock \emph{Journal of Computational Physics}, 443:\penalty0 110525, 2021.

\bibitem[Fasel et~al.(2022)Fasel, Kutz, Brunton, and
  Brunton]{fasel2022ensemble}
Urban Fasel, J.~Nathan Kutz, Bingni~W. Brunton, and Steven~L. Brunton.
\newblock Ensemble-{SINDy}: Robust sparse model discovery in the low-data,
  high-noise limit, with active learning and control.
\newblock \emph{Proceedings of the Royal Society A}, 478\penalty0
  (2260):\penalty0 20210904, 2022.

\bibitem[Kamienny et~al.(2022)Kamienny, d'Ascoli, Lample, and
  Charton]{kamienny2022end}
Pierre-Alexandre Kamienny, St{\'e}phane d'Ascoli, Guillaume Lample, and
  Fran{\c{c}}ois Charton.
\newblock End-to-end symbolic regression with transformers.
\newblock In \emph{Advances in Neural Information Processing Systems},
  volume~35, 2022.

\bibitem[d'Ascoli et~al.(2024)d'Ascoli, Becker, Mathis, Schwaller, and
  Kilbertus]{dascoli2024odeformer}
St{\'e}phane d'Ascoli, S{\"o}ren Becker, Alexander Mathis, Philippe Schwaller,
  and Niki Kilbertus.
\newblock {ODEFormer}: Symbolic regression of dynamical systems with
  transformers.
\newblock In \emph{International Conference on Learning Representations
  (ICLR)}, 2024.

\bibitem[Sun et~al.(2024)Sun, Liu, Zhang, and Schaeffer]{sun2024prosepde}
Jingmin Sun, Yuxuan Liu, Zecheng Zhang, and Hayden Schaeffer.
\newblock Towards a foundation model for partial differential equations:
  Multi-operator learning and extrapolation.
\newblock \emph{arXiv preprint arXiv:2404.12355}, 2024.

\bibitem[Liu et~al.(2026)Liu, Luo, Zhu, Yang, and Lu]{liu2026foundpde}
Langchen Liu, Junlin Luo, Min Zhu, Zhuoran Yang, and Lu~Lu.
\newblock {FoundPDE}: A generative foundation model for data-driven symbolic
  discovery of partial differential equations.
\newblock \emph{arXiv preprint arXiv:2610.26963}, 2026.

\bibitem[Xu et~al.(2025)Xu, Chen, Cao, Tang, Du, Li, Callaghan, and
  Zhang]{xu2025eqgpt}
Hao Xu, Yuntian Chen, Rui Cao, Tianning Tang, Mengge Du, Jian Li, Adrian~H.
  Callaghan, and Dongxiao Zhang.
\newblock Generative discovery of partial differential equations by learning
  from math handbooks.
\newblock \emph{Nature Communications}, 16:\penalty0 10255, 2025.

\bibitem[Shojaee et~al.(2025)Shojaee, Meidani, Gupta, Farimani, and
  Reddy]{shojaee2025llmsr}
Parshin Shojaee, Kazem Meidani, Shashank Gupta, Amir~Barati Farimani, and
  Chandan~K. Reddy.
\newblock {LLM-SR}: Scientific equation discovery via programming with large
  language models.
\newblock In \emph{International Conference on Learning Representations
  (ICLR)}, 2025.

\bibitem[Both et~al.(2021)Both, Choudhury, Sens, and Kusters]{both2021deepmod}
Gert-Jan Both, Subham Choudhury, Pierre Sens, and Remy Kusters.
\newblock {DeepMoD}: Deep learning for model discovery in noisy data.
\newblock \emph{Journal of Computational Physics}, 428:\penalty0 109985, 2021.

\bibitem[Li et~al.(2026)Li, Cui, Qi, Zhang, Li, and Yin]{weakpdenet2026}
Xinxin Li, Xingyu Cui, Jin Qi, Juan Zhang, Da~Li, and Junping Yin.
\newblock {Weak-PDE-Net}: Discovering open-form {PDE}s via differentiable
  symbolic networks and weak formulation.
\newblock \emph{arXiv preprint arXiv:2603.22951}, 2026.

\bibitem[Bideh et~al.(2026)Bideh, Georgievska, and Gryak]{bideh2026mdbench}
Amirmohammad~Ziaei Bideh, Aleksandra Georgievska, and Jonathan Gryak.
\newblock {MDBench}: Benchmarking data-driven methods for model discovery.
\newblock In \emph{Proceedings of the AAAI Conference on Artificial
  Intelligence}, volume~40, pages 19746--19754, 2026.

\bibitem[Krull et~al.(2019)Krull, Buchholz, and Jug]{krull2019noise2void}
Alexander Krull, Tim-Oliver Buchholz, and Florian Jug.
\newblock {Noise2Void}: Learning denoising from single noisy images.
\newblock In \emph{Proceedings of the IEEE/CVF Conference on Computer Vision
  and Pattern Recognition}, pages 2129--2137, 2019.

\bibitem[Koehler et~al.(2024)Koehler, Niedermayr, Westermann, and
  Thuerey]{apebench}
Felix Koehler, Simon Niedermayr, R{\"u}diger Westermann, and Nils Thuerey.
\newblock {APEBench}: A benchmark for autoregressive neural emulators of
  {PDE}s.
\newblock In \emph{Advances in Neural Information Processing Systems (Datasets
  and Benchmarks Track)}, volume~37, 2024.

\bibitem[Reinbold et~al.(2020)Reinbold, Gurevich, and
  Grigoriev]{reinbold2020using}
Patrick A.~K. Reinbold, Daniel~R. Gurevich, and Roman~O. Grigoriev.
\newblock Using noisy or incomplete data to discover models of spatiotemporal
  dynamics.
\newblock \emph{Physical Review E}, 101:\penalty0 010203, 2020.

\end{thebibliography}

\clearpage
\appendix

\section{Data Generation Details}
\label{app:data}

Table~\ref{tab:sim-settings} lists the simulation settings of each family. Every family uses $N_x = 128$ grid points on a periodic domain and stores 501 frames (500 steps). Linear families are advanced exactly in Fourier space; Burgers, KdV, and KS use second-order exponential time differencing (ETDRK2) with $2/3$ dealiasing, as implemented in APEBench. When the stride exceeds one, the solver runs at the finer step and every stride-th frame is stored. Each initial condition is a random truncated Fourier series
\[
u_0(x) = \sum_{k=1}^{5} a_k \sin\!\Bigl(\frac{2\pi k x}{L} + \phi_k\Bigr),
\qquad a_k \sim U(-1,1),\quad \phi_k \sim U(0, 2\pi),
\]
with zero mean, rescaled so that $\max_x |u_0| = 1$, and drawn with its own random seed. For KS, the solver first runs the stated number of warm-up steps, and the stored trajectory starts afterwards.

\begin{table}[!htbp]
\centering
\small
\setlength{\tabcolsep}{5pt}
\caption{Simulation settings per family: domain length $L$, solver step, storage stride, stored time step $\Delta t$, duration $T_{\mathrm{total}}$, solver warm-up steps before the first stored frame, and time integrator.}
\label{tab:sim-settings}
\begin{tabular}{lrrrrrrl}
\toprule
Family & $L$ & Solver step & Stride & $\Delta t$ & $T_{\mathrm{total}}$ & Warm-up & Integrator \\
\midrule
Advection & 1 & 0.01 & 1 & 0.01 & 5.0 & 0 & Exact \\
Diffusion & 1 & 0.004 & 1 & 0.004 & 2.0 & 0 & Exact \\
Advection--diffusion & 1 & 0.01 & 1 & 0.01 & 5.0 & 0 & Exact \\
Burgers & 1 & 0.01 & 1 & 0.01 & 5.0 & 0 & ETDRK2 \\
KdV & 50 & 0.0005 & 12 & 0.006 & 3.0 & 0 & ETDRK2 \\
KS & 60 & 0.01 & 4 & 0.04 & 20.0 & 2000 & ETDRK2 \\
\bottomrule
\end{tabular}
\end{table}

\section{Model and Training Details}
\label{app:model}

\paragraph{Terminology.}
Table~\ref{tab:terminology} fixes the three units of data used throughout the paper. Identification is always per trajectory: every method, including each baseline, receives one noisy trajectory and returns one equation.

\begin{table}[h]
\centering
\caption{Terminology. Each level contains many instances of the level above.}
\label{tab:terminology}
\small
\begin{tabularx}{\linewidth}{@{}>{\raggedright\arraybackslash}p{0.17\linewidth} >{\raggedright\arraybackslash}X >{\raggedright\arraybackslash}p{0.17\linewidth} >{\raggedright\arraybackslash}p{0.27\linewidth}@{}}
\toprule
Term & Definition & Fixed & Example \\
\midrule
Trajectory & one simulated space--time field $u(x,t)$ & terms, coefficients, and IC & one field in Figure~\ref{fig:pde-gallery} \\
System & one governing equation, observed through trajectories that differ in their IC & terms and coefficients & $u_t = -0.4\,u u_x + 0.003\,u_{xx}$ \\
Family (PDE class) & an equation form with free coefficients & terms only & Burgers: $u_t = -b\,u u_x + \nu\,u_{xx}$ \\
\bottomrule
\end{tabularx}
\end{table}

\paragraph{Hardware and software.}
All experiments---data generation, denoiser and encoder training, inference, and the timing measurements in Table~\ref{tab:efficiency}---were run on a single Apple M3 Pro laptop with 18\,GB of unified memory under macOS 26.6.2, using the Metal (MPS) backend for the neural components. The software is Python 3.12 with PyTorch 2.12.0, JAX 0.10.1, APEBench 0.1.1 (Exponax 0.1.0), PySINDy 2.1.0, and PySR 1.5.10 (Julia 1.13.2, SymbolicRegression.jl 1.11.3).

\paragraph{Scope.}
The noisy pipeline itself (denoiser, encoder, weak-form solve) is specified in the main text. This appendix records the input representation and encoder architecture shared by all configurations, the individual loss terms, and the clean-data reference path, labeled UNIDE (clean specialist) in the clean reference table (Appendix~\ref{app:clean}), where coefficients are identified by pointwise full-window least squares instead of the weak form. All methods are evaluated per trajectory on the same stored fields: strong-form derivative features and pointwise least-squares refits use the 499 interior frames left by the centered time difference, the weak-form solves integrate over the full stored field, and the encoder reads a centered 64-frame crop (input representation below).

\paragraph{Denoiser architecture.}
The denoiser $D_\eta$ is a residual 2D CNN over $(t,x)$: a $1\times1$ input convolution to width 48, eight $3\times3$ residual blocks with dilations cycling $(1,2,4,8)$, and a $1\times1$ output convolution. Padding is circular in $x$ (periodic domain) and replicate in $t$. Training hides a random $2\%$ of space--time pixels per step (the blind spot): each masked entry is replaced in the input by the value of the pixel offset by a random one to three grid points in both $t$ and $x$ (one offset drawn per training step), and the loss is evaluated only at the masked entries against the original noisy values, on random 96-frame time crops. Optimization uses Adam with learning rate $5\times10^{-4}$, cosine annealing, gradient clipping at $1.0$, and batch size 16. Each training crop receives a fresh noise realization at a level drawn from $U(0,\sigma)$ at every step, so the denoiser sees many independent noisy copies of each simulated training trajectory but never a clean target. Each training crop is standardized by its own mean and standard deviation. At inference, the full noisy field is standardized in the same way and processed in 256-frame windows that start every 240 frames, with a final window aligned to the end of the field (for 501 frames, windows start at frames 0, 240, and 245); predictions are averaged where windows overlap and mapped back to physical units.

\begin{figure}[!htbp]
\centering
\includegraphics[width=0.7\linewidth]{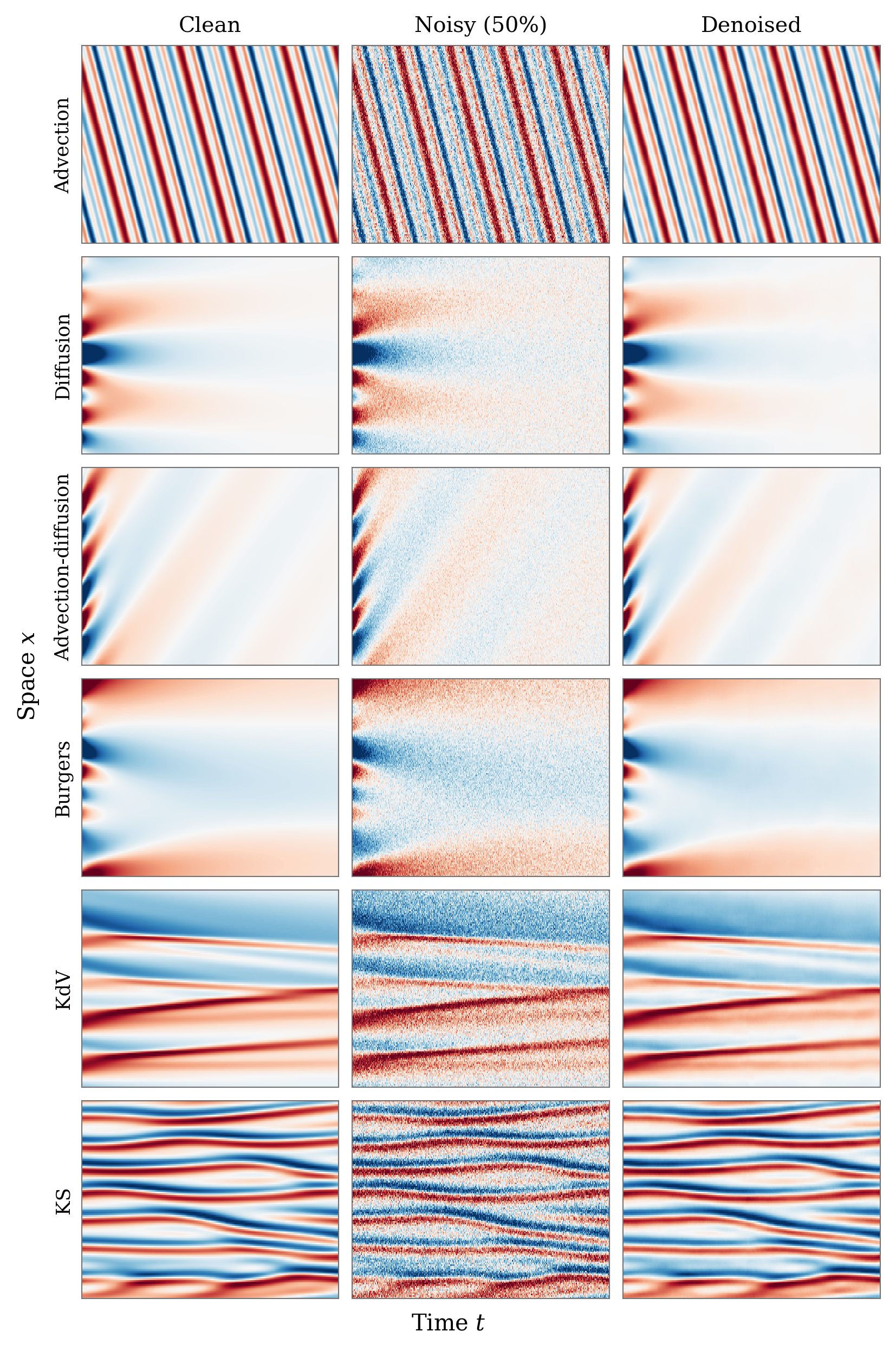}
\caption{Self-supervised denoising across the six PDE families at $50\%$ observation noise. Each row is one held-out test trajectory (horizontal: time, vertical: space); columns show the clean field, the noisy observation, and the denoiser output. The color scale is symmetric and set per row from the clean field.}
\label{fig:denoise-fields}
\end{figure}

\begin{figure}[!htbp]
\centering
\includegraphics[width=0.68\linewidth]{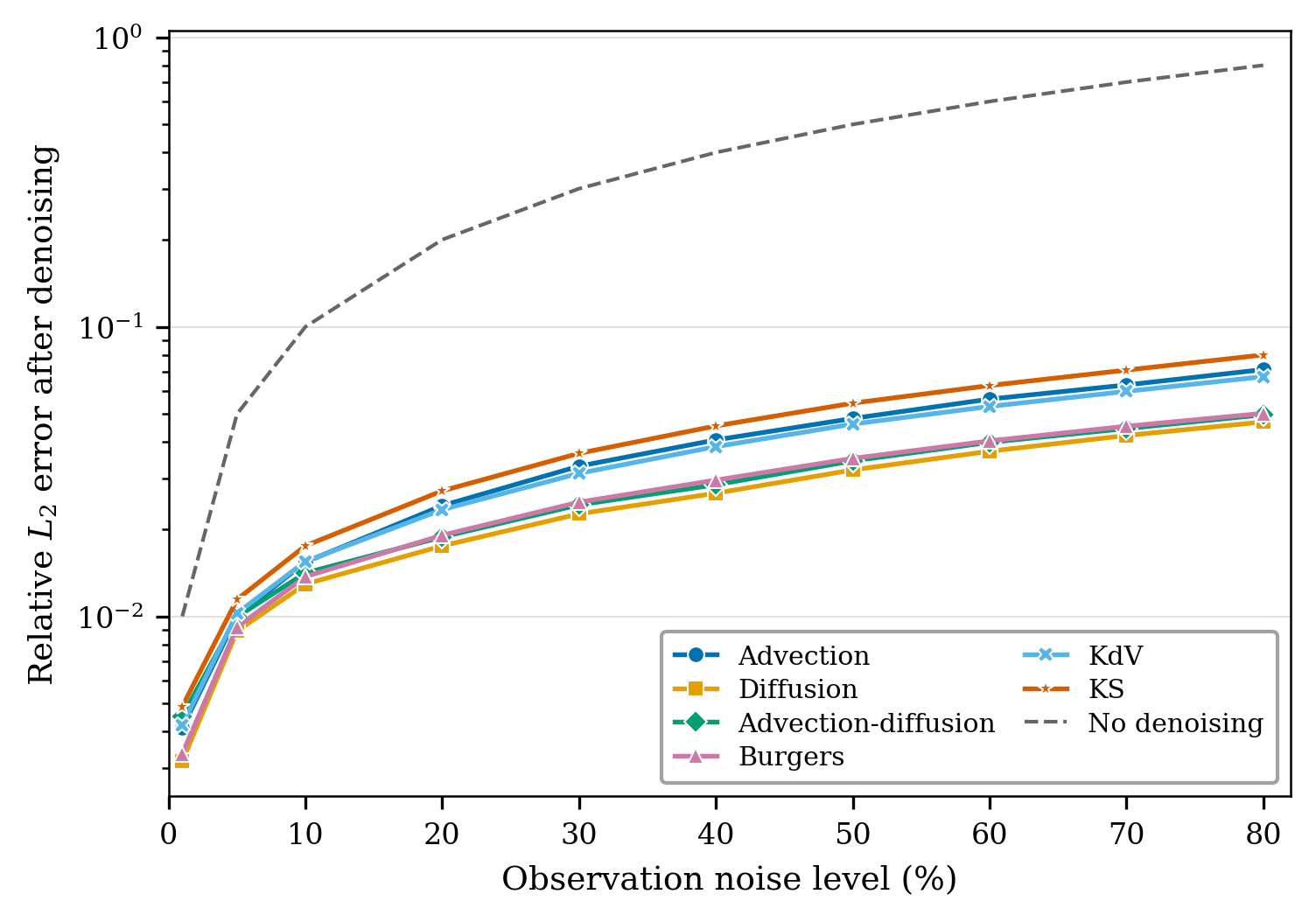}
\caption{Residual relative $L_2$ error after denoising versus observation-noise level, one line per PDE family (log scale); each point is the mean over the family's 150 held-out test trajectories. Each level uses its own level-matched specialist, and the dashed line is the no-denoising identity (residual error $\approx$ noise level).}
\label{fig:denoise-summary}
\end{figure}

\paragraph{Denoising quality.}
Figure~\ref{fig:denoise-fields} shows the denoiser on one held-out test trajectory per family at $50\%$ observation noise: it recovers each family's space--time structure---traveling fronts, dispersive trains, chaotic texture---from the noisy observations, alongside the clean reference fields. Figure~\ref{fig:denoise-summary} quantifies this across the full sweep, averaged over the 150 held-out test trajectories of each family (values summarized in Section~\ref{subsec:denoise}). From $10\%$ noise onward the dissipative families (diffusion, advection--diffusion, Burgers) clean up best, while advection, KdV, and KS, whose fine-scale content persists undamped, retain the most residual---the residual the weak-form coefficient solve is built to tolerate.

\paragraph{Input representation.}
For one sample, the model receives
\[
X \in \mathbb{R}^{T \times 12 \times N_x}, \qquad
g \in \mathbb{R}^{T \times 3}, \qquad
z \in \mathbb{R}^{d_z},
\]
where $T=64$ is the temporal crop length and $N_x=128$ is the spatial resolution. The twelve feature channels are
\[
u_t,\;
u,\;
u^2,\;
u^3,\;
u_x,\;
u_{xx},\;
u_{xxx},\;
u_{xxxx},\;
u u_x,\;
u^2u_x,\;
u\,u_{xx},\;
u_x^2.
\]
In the clean configuration, the first channel is the centered finite-difference time derivative computed from the stored trajectory and the remaining channels come from the precomputed library tensor; in the noisy configuration, all twelve channels are recomputed from the denoised field $\tilde u$. The time features are
\[
g_t = [t,\; t/T_{\mathrm{total}},\; \Delta t],
\]
and the metadata vector is
\[
z = [\Delta x,\; \Delta t,\; T_{\mathrm{total}},\; L,\; N_x,\; T_{\mathrm{obs}}].
\]
Here $T_{\mathrm{obs}}$ is the number of observed interior frames, so $d_z=6$; every encoder in this paper, clean and noisy, uses this metadata vector.

\paragraph{Trajectory encoder.}
The per-frame spatial CNN $f_{\mathrm{cnn}}$, the time embedding $f_{\mathrm{time}}$, and the transformer $f_{\mathrm{tr}}$ are defined in Section~\ref{subsec:support}. Here $f_{\mathrm{cnn}}$ is two 1D convolution layers with GELU activations, global average pooling over space, a linear projection, and LayerNorm, mapping each frame $X_t \in \mathbb{R}^{12 \times N_x}$ to a code $s_t \in \mathbb{R}^{d}$; a sinusoidal positional encoding is added to $r_{1:T}$ before the transformer. The trajectory representation is obtained by attention pooling over the frames,
\[
\alpha_t = \frac{\exp(w^\top h_t)}{\sum_{s=1}^T \exp(w^\top h_s)},
\qquad
h_{\mathrm{traj}} = \sum_{t=1}^T \alpha_t h_t.
\]
Metadata is embedded separately as $h_{\mathrm{meta}} = f_{\mathrm{meta}}(z)$, and the final shared representation is
\[
h = \operatorname{GELU}\!\left(W [h_{\mathrm{traj}}; h_{\mathrm{meta}}] + b\right).
\]
The spatial CNN uses 64 and 128 channels with kernel sizes 7 and 5; the time and metadata MLPs use one hidden layer of 64 and 32 units; and the transformer uses $d_{\mathrm{model}}=128$, 4 heads, 2 layers, and dropout 0.1, for about 0.5M parameters in total. The encoder is trained with AdamW (learning rate $10^{-3}$ halved every 20 epochs, weight decay $10^{-4}$) and batch size 32; the noisy-configuration encoders are warm-started from the clean-trained model.

\paragraph{Auxiliary coefficient head and clean-reference solve.}
Let $\Theta \in \mathbb{R}^{(T N_x)\times M}$ be the reshaped library matrix for the crop, let $b \in \mathbb{R}^{T N_x}$ denote the reshaped $u_t$, and let $r \in \mathbb{R}^{M}$ be the per-column training-set RMS statistics; write $\Theta_n = \Theta\operatorname{Diag}(r)^{-1}$ for the column-normalized library. During training, an auxiliary linear head maps $h$ to raw coefficients $\tilde{\xi}$, which are gated by the (detached) term probabilities and unscaled to physical units:
\[
\hat{\xi}_n = \sigma(\ell) \odot \tilde{\xi}, \qquad
\hat{\xi}_{\mathrm{direct}} = \operatorname{Diag}(r)^{-1}\hat{\xi}_n.
\]
The gating forces inactive terms toward zero. These auxiliary coefficients exist only to drive $\mathcal{L}_{\mathrm{coef}}$, $\mathcal{L}_{\mathrm{res}}$, and $\mathcal{L}_{\mathrm{sp}}$ during training; no reported coefficient comes from this head.

At evaluation time, the probabilities are thresholded to a hard mask
\[
\hat{m}_j = \mathbf{1}\{\sigma(\ell_j) > \tau\},
\]
where $\tau$ is selected on validation by exact match over $\{0.3, 0.4, \ldots, 0.9\}$ with ties broken toward $0.5$; the selected value is $0.5$ in both clean settings (the noisy pipeline tunes $\tau$ by F1 at each noise level, Section~\ref{subsec:support}). The clean-reference coefficients are then computed by
\[
\hat{\xi}
=
\arg\min_{\xi}
\left\|
\Theta \operatorname{Diag}(\hat{m})\xi - b
\right\|_2^2,
\]
where $\Theta$ and $b$ here are the full-window library and $u_t$ over all 499 interior frames, not the crop-based versions that drive the training losses. The system is solved with a float64 SVD-based least-squares routine.

\paragraph{Loss terms.}
Let $m \in \{0,1\}^{M}$ and $\xi \in \mathbb{R}^{M}$ denote the true support and coefficients. The individual terms of the objective in Section~\ref{subsec:objective} are defined as follows. The term loss is weighted binary cross-entropy on the logits:
\[
\mathcal{L}_{\mathrm{term}}
=
\frac{1}{M}\sum_{j=1}^{M}
\left[
-w^+ m_j \log \sigma(\ell_j)
-
(1-m_j)\log\!\bigl(1-\sigma(\ell_j)\bigr)
\right].
\]
The coefficient term is relative MAE on active terms only:
\[
\mathcal{L}_{\mathrm{coef}}
=
\frac{1}{\max(1,\|m\|_0)}
\sum_{j=1}^{M}
m_j
\frac{|\hat{\xi}_{\mathrm{direct},j} - \xi_j|}{\max(|\xi_j|, 10^{-2})}.
\]
The residual term compares the predicted right-hand side to the observed $u_t$, with both norms averaged over the $N$ space--time entries:
\[
\mathcal{L}_{\mathrm{res}}
=
\frac{
\tfrac{1}{N}\left\|
\Theta_n \hat{\xi}_n - b
\right\|_2^2
}{
\max\bigl(\tfrac{1}{N}\|b\|_2^2,\,10^{-8}\bigr)
}.
\]
The sparsity term is
\[
\mathcal{L}_{\mathrm{sp}} = \|\hat{\xi}_{\mathrm{direct}}\|_1 / M,
\]
and $\mathcal{L}_{\mathrm{fam}}$ is the auxiliary cross-entropy loss for PDE-family classification. The term loss uses positive-class weighting with $w^+=5$. In the implementation, the normalizations of $\mathcal{L}_{\mathrm{coef}}$ and $\mathcal{L}_{\mathrm{res}}$ pool over the training batch rather than per sample. In the noisy configuration the denoiser trains for 120 epochs and the encoder for 60; the clean-reference configuration trains for 100 epochs. For the denoiser, the checkpoint with the lowest validation loss is used. For the encoder, the final epoch is used if its validation Term F1 is within 0.01 of the best epoch's, and otherwise the best-F1 epoch, with ties broken by lower validation coefficient error.

\paragraph{Weak-form least squares: full formulation.}
\label{app:weakform}
The weak library represents ten of the twelve library terms: seven of the form $\partial_x^{d} f(u)$ with state polynomials $f(u) \in \{u, u^2, u^3\}$ and derivative orders $d \le 4$, and the three interaction terms treated below. The remaining two ($1$ and $u_x^2$) appear in no true support; the weak library omits them, so they are excluded from the pruned support and receive coefficient zero (in practice the encoder never selects them on any test trajectory). Following the weak formulation of WSINDy~\citep{messenger2021weak}, the dynamics are tested against $K$ compactly supported functions on local space--time subdomains, placed at random as in \citet{reinbold2020using}; we use PySINDy's \texttt{WeakPDELibrary}. Subdomain $k$ is a rectangle $\Omega_k$ whose center $(x_k, t_k)$ is drawn uniformly at random over the interior of the space--time domain (fixed seed), with nominal half-widths $H_x$ and $H_t$ of one twentieth of the grid span in each direction; the rectangle is then shrunk so that its boundaries lie on grid points, and its center $(x_k, t_k)$ and half-widths $H_{x,k}$ and $H_{t,k}$ are recomputed from the shrunk rectangle. Its test function is the separable polynomial bump
\[
\psi_k(x,t) = \bigl(1-\tilde{x}^2\bigr)^{p} \bigl(1-\tilde{t}^2\bigr)^{p},
\qquad
\tilde{x} = \frac{x - x_k}{H_{x,k}},
\quad
\tilde{t} = \frac{t - t_k}{H_{t,k}},
\qquad p = 4,
\]
which vanishes to order $p$ on $\partial\Omega_k$, so every boundary term through fourth-order derivatives vanishes under integration by parts; $p=4$ is the minimal such degree. The fixed $H$ and $p$ are the implementation defaults---WSINDy itself selects both adaptively. Multiplying $u_t = \sum_j \xi_j\,\partial_x^{d_j} f_j(u)$ by $\psi_k$, integrating over $\Omega_k$, and moving all derivatives onto $\psi_k$ gives row $k$ of the weak system $\Theta_{\mathrm{w}}\xi \approx b_{\mathrm{w}}$:
\[
b_{\mathrm{w},k} = -\int_{\Omega_k} (\partial_t \psi_k)\, \tilde{u}\; dx\,dt,
\qquad
[\Theta_{\mathrm{w}}]_{kj} = (-1)^{d_j} \int_{\Omega_k} \bigl(\partial_x^{d_j} \psi_k\bigr)\, f_j(\tilde{u})\; dx\,dt .
\]
The three interaction columns use a partially integrated form that keeps a first-order finite difference of the observed field under the integral: $\int_{\Omega_k}\psi_k f(\tilde u)\partial_x^{\mathrm{FD}}\tilde u$ for $u u_x$ and $u^2 u_x$. Because $u u_x=\tfrac12\partial_x(u^2)$ and $u^2u_x=\tfrac13\partial_x(u^3)$, these two columns also admit the fully weak integrals $-\tfrac12\int_{\Omega_k}(\partial_x\psi_k)\tilde u^2$ and $-\tfrac13\int_{\Omega_k}(\partial_x\psi_k)\tilde u^3$; all reported experiments use the partially integrated form. For $u u_{xx}$, the integral is $-\int_{\Omega_k}[\psi_k(\partial_x^{\mathrm{FD}}\tilde u)^2+(\partial_x\psi_k)\tilde u\partial_x^{\mathrm{FD}}\tilde u]$. Pure derivative columns move all derivatives onto the test function. The integrals use the library's analytic weights for piecewise-linear interpolation of the grid values. Selected columns are normalized to unit root-mean-square, fitted by least squares, and rescaled to the original units; inactive terms receive zero coefficients. UNIDE uses $K=400$ subdomains. The WSINDy baseline selects its support with $K=800$ before its $K=400$ coefficient refit.

\paragraph{Residual pruning inside the solve.}
Let $r(S) = \|b_{\mathrm{w}} - \Theta_{\mathrm{w},S}\,\hat\xi_S\|_2 / \|b_{\mathrm{w}}\|_2$ denote the normalized residual of the weak-form fit on an active set $S$, computed with the same column normalization as the solve. Starting from $S = \{j : \hat m_j = 1\}$, the solve repeatedly finds the term whose single-term removal increases the residual least and removes it while
\[
r\bigl(S \setminus \{j\}\bigr) \le (1+\delta)\, r(S),
\]
stopping as soon as no removal passes the test (at least one term is always kept). The tolerance $\delta$ is tuned per noise level on the validation split by exact match (ties broken toward lower coefficient error), over the grid $\{0, 0.004, 0.008, 0.012, 0.02, 0.05\}$, where $\delta = 0$ disables the pruning, so levels that do not benefit keep the unpruned solve automatically. The two tuned quantities deliberately use different criteria: the mask threshold $\tau_{\mathrm{val}}$ is tuned by F1, which favors recall and lets the encoder over-include, while $\delta$ is tuned by exact match, so the selected tolerance trades support precision against recall on validation data.

\section{Evaluation Details}
\label{app:eval}

\subsection{Metric Definitions}

For sample $i$, let $m_i \in \{0,1\}^M$ be the true support, $\hat m_i \in \{0,1\}^M$ the predicted support, $A_i = \{j : m_{ij}=1\}$ the true active set, and $\xi_i, \hat\xi_i \in \mathbb{R}^M$ the true and predicted coefficients. All metrics are computed per sample and then averaged over the evaluation set.

\paragraph{Structure metrics.}
Define
\[
\mathrm{TP}_i = \sum_{j=1}^M \mathbf{1}\{m_{ij}=1 \wedge \hat m_{ij}=1\},\quad
\mathrm{FP}_i = \sum_{j=1}^M \mathbf{1}\{m_{ij}=0 \wedge \hat m_{ij}=1\},
\]
\[
\mathrm{FN}_i = \sum_{j=1}^M \mathbf{1}\{m_{ij}=1 \wedge \hat m_{ij}=0\}.
\]
The per-sample precision, recall, and F1 are
\[
P_i = \frac{\mathrm{TP}_i}{\mathrm{TP}_i+\mathrm{FP}_i+\varepsilon},\qquad
R_i = \frac{\mathrm{TP}_i}{\mathrm{TP}_i+\mathrm{FN}_i+\varepsilon},
\qquad
\mathrm{F1}_i = \frac{2P_iR_i}{P_i+R_i+\varepsilon},
\]
with a small $\varepsilon$ for numerical stability. The reported \textbf{Term F1} is the mean of $\mathrm{F1}_i$ over samples, and \textbf{Exact match} is
\[
\mathrm{Exact} = \frac{1}{N}\sum_{i=1}^N \mathbf{1}\{\hat m_i = m_i\},
\]
the fraction of samples whose support is identified perfectly.

\paragraph{Coefficient metrics.}
Coefficient errors are evaluated only on the true active terms:
\[
\mathrm{RMSE}_i =
\sqrt{
\frac{1}{|A_i|}
\sum_{j \in A_i} (\hat\xi_{ij}-\xi_{ij})^2
},
\qquad
\mathrm{RelErr}_i =
\frac{1}{|A_i|}
\sum_{j \in A_i}
\frac{|\hat\xi_{ij}-\xi_{ij}|}{\max(|\xi_{ij}|,\varepsilon)},
\]
reported as \textbf{Coefficient RMSE} and \textbf{coefficient relative error}, with $\varepsilon=10^{-8}$ used only for numerical stability. The noisy residual-pruning evaluations, baseline sweeps, PySR runs, and matched controls use the denominator $\max(|\xi_{ij}|,\varepsilon)$; other neural evaluations use $|\xi_{ij}|+\varepsilon$. The two conventions differ negligibly at the reported precision, because every active coefficient in the benchmark has magnitude at least $10^{-3}$.

\paragraph{Reading the two metric families together: a worked example.}
Because coefficient error is computed only on the true active terms, a spurious extra term hurts the structure metrics but is invisible to coefficient error. Advection illustrates this. PySR can occasionally improve the fitted coefficient of the true transport term by adding a spurious high-order derivative that absorbs part of the finite-difference $u_t$ bias; coefficient error then stays low even though exact match drops. Conversely, when a method misses a true active term, both structure and coefficient metrics degrade together. This is why we treat support identification as the primary metric and read coefficient quality in that context.

\subsection{Baseline Implementations}
\label{app:baselines}

All three baseline implementations receive the same stored $128\times501$ trajectories, with features recomputed from noisy observations. Their internal search spaces differ. The pointwise PySINDy PDE library contains 20 columns and projects its fitted coefficients onto 11 of the 12 library terms. The weak library contains 19 columns and projects onto 10 output terms. Each mapping discards nine additional interaction columns; these columns can affect selection before projection, and an omitted selected column is not counted as a false positive in the reported metrics. UNIDE predicts 12 term probabilities and its weak solve retains the ten terms of the weak library. The PySR feature-selection baseline receives all 12 pointwise library columns. These results are therefore comparisons of the stated implementations, not an exactly matched candidate-library experiment. Sparse thresholds are tuned by validation F1 using 48 trajectories (8 per family), whereas UNIDE's mask threshold uses the full 300-trajectory validation split; the matched controls of Section~\ref{sec:matched-controls} tune weak-form STLSQ on the full split over an eleven-value threshold grid and still trail UNIDE in exact match. PySR's search settings are fixed.

\paragraph{SINDy.}
Implemented with PySINDy: a pointwise PDE library with a constant column, polynomial state functions up to degree~3, periodic finite-difference spatial derivatives up to order~4, and their interaction terms (20 columns), and the STLSQ optimizer (ridge parameter $\alpha=10^{-5}$, normalized columns). The threshold is tuned over $\{0.5, 1, 2, 5\}$. Coefficients on the selected support are de-biased by an ordinary least-squares refit on the pointwise library columns.

\paragraph{WSINDy.}
Uses the weak system of Appendix~\ref{app:weakform} built with $K=800$ subdomains, and STLSQ in weak form. The threshold is tuned over $\{0.005, 0.008, 0.01\}$. Coefficients are then re-solved by weak-form least squares on the method's own selected support ($K=400$, the same least-squares solve as UNIDE but without the residual pruning, which is part of UNIDE's method rather than of the evaluation), because the standard pointwise refit collapses under noise and would understate the baseline.

\paragraph{PySR.}
Per-trajectory symbolic regression over the 12 library columns computed from the observed field as input features, with target $u_t$: binary operators $\{+,\times\}$, no unary operators, 40 iterations, maximum equation size 25, parsimony $3.2\times10^{-3}$, and \texttt{model\_selection=score}, run deterministically with a fixed seed. Feature columns are standardized to unit RMS for the search, and $4{,}000$ space--time samples are drawn per trajectory for tractability. Support is defined as the set of input library variables appearing in the selected expression. We then discard the expression's nonlinear composition and refit a linear combination of those variables by least squares on all observed library rows. PySR was run on every test trajectory (10,620: the random-trajectory test set at the eleven stated noise levels and the coefficient-extrapolation test set at 0 and 10\%), with no fit failures; its input library is recomputed from the observed field at every level, including zero noise. Thus this is a PySR-based feature-selection baseline, not an evaluation of unrestricted symbolic-expression recovery. In particular, multiplying two library variables in the search does not add their product as a new scored output term.

\paragraph{Matched controls.}
All weak systems in the controls of Section~\ref{sec:matched-controls} use the same ten candidate terms and 400 windows, with selected terms scored in the full twelve-term library. STLSQ thresholds are selected on the full validation partition from $\{0.001,\allowbreak 0.003,\allowbreak 0.005,\allowbreak 0.008,\allowbreak 0.01,\allowbreak 0.02,\allowbreak 0.05,\allowbreak 0.1,\allowbreak 0.2,\allowbreak 0.5,\allowbreak 1\}$, with ridge parameter $10^{-5}$ and column normalization. The joint STLSQ-plus-pruning variant searches that threshold grid and the existing pruning grid by validation F1; selecting instead by validation exact match, the criterion used for UNIDE's pruning tolerance, picks the same configuration in all three settings.

\subsection{Clean-Data Reference Results}
\label{app:clean}

Table~\ref{tab:clean-ref} reports the 0\% noise reference for both settings. On clean data the per-trajectory baselines are strong, PySR most of all (exact match 0.953). UNIDE's clean specialist recovers structure more accurately still in both settings, and on the random-trajectory set its coefficients as well; under clean coefficient extrapolation PySR fits the coefficients best. \textbf{UNIDE (clean specialist)} denotes the configuration UNIDE uses when there is no noise to remove: the denoiser is skipped, the shared encoder trained on clean features predicts the support, and coefficients come from pointwise full-window least squares on that support, with no weak form and hence no residual pruning (Appendix~\ref{app:model}). This replaces the denoise--encode--weak-form pipeline used at every noisy level, because there is no noise to remove and the weak form is needed only to resist noise.

\begin{table}[!htbp]
\centering
\small
\setlength{\tabcolsep}{4pt}
\caption{Clean-data (0\% noise) reference results. Top: the random-trajectory test set (900 trajectories). Bottom: clean coefficient extrapolation (360 trajectories, held-out settings with the largest coefficient-vector norms). UNIDE (clean specialist) is defined in the text above. Best per metric in bold.}
\label{tab:clean-ref}
\begin{tabular}{lcccc}
\toprule
Method & Exact $\uparrow$ & Term F1 $\uparrow$ & Coef. RMSE $\downarrow$ & Coef. rel\% $\downarrow$ \\
\midrule
\multicolumn{5}{l}{\emph{Random-trajectory test set (clean)}} \\
SINDy & 0.690 & 0.884 & 1.31e-02 & 1.24 \\
WSINDy & 0.836 & 0.938 & 2.70e-02 & 2.29 \\
PySR & 0.953 & 0.984 & 8.45e-03 & 0.50 \\
\textbf{UNIDE (clean specialist)} & \textbf{0.999} & \textbf{1.000} & \textbf{1.03e-03} & \textbf{0.21} \\
\midrule
\multicolumn{5}{l}{\emph{Coefficient extrapolation (clean)}} \\
SINDy & 0.503 & 0.809 & 1.12e-02 & 1.06 \\
WSINDy & 0.778 & 0.918 & 4.74e-02 & 2.25 \\
PySR & 0.931 & 0.977 & \textbf{1.83e-03} & \textbf{0.41} \\
\textbf{UNIDE (clean specialist)} & \textbf{0.989} & \textbf{0.994} & 4.97e-03 & 0.73 \\
\bottomrule
\end{tabular}
\end{table}

\section{Noise-Sweep Results}
\label{app:sweep}

Table~\ref{tab:noise-sweep} lists the values plotted in Figure~\ref{fig:main-results} for every method and noisy level.

\begin{table}[!htbp]
\centering
\small
\setlength{\tabcolsep}{4pt}
\caption{Noise sweep on the random-trajectory test set (900 trajectories per row); these are the values plotted in Figure~\ref{fig:main-results}. UNIDE uses the specialist trained at each level, with the pruned solve; the baselines receive the raw noisy fields. Best per metric within each noise level in bold.}
\label{tab:noise-sweep}
\begin{tabular}{llcccc}
\toprule
Noise & Method & Exact $\uparrow$ & Term F1 $\uparrow$ & Coef. RMSE $\downarrow$ & Coef. rel\% $\downarrow$ \\
\midrule
1\% & SINDy & 0.016 & 0.471 & 5.00e-01 & 50.64 \\
 & WSINDy & 0.820 & \textbf{0.941} & \textbf{2.79e-02} & \textbf{2.84} \\
 & PySR & 0.167 & 0.517 & 6.51e-01 & 64.57 \\
 & UNIDE & \textbf{0.866} & 0.938 & 5.15e-02 & 10.32 \\
\midrule
5\% & SINDy & 0.000 & 0.319 & 5.89e-01 & 70.41 \\
 & WSINDy & 0.566 & 0.789 & \textbf{2.56e-02} & \textbf{3.59} \\
 & PySR & 0.150 & 0.417 & 7.18e-01 & 80.69 \\
 & UNIDE & \textbf{0.877} & \textbf{0.942} & 4.27e-02 & 10.18 \\
\midrule
10\% & SINDy & 0.000 & 0.307 & 6.70e-01 & 77.51 \\
 & WSINDy & 0.517 & 0.746 & 2.99e-02 & 6.19 \\
 & PySR & 0.096 & 0.334 & 7.62e-01 & 86.55 \\
 & UNIDE & \textbf{0.960} & \textbf{0.984} & \textbf{2.92e-02} & \textbf{4.99} \\
\midrule
20\% & SINDy & 0.003 & 0.318 & 7.66e-01 & 86.86 \\
 & WSINDy & 0.429 & 0.706 & 4.45e-02 & 13.52 \\
 & PySR & 0.007 & 0.255 & 8.08e-01 & 90.24 \\
 & UNIDE & \textbf{0.971} & \textbf{0.989} & \textbf{3.07e-02} & \textbf{4.46} \\
\midrule
30\% & SINDy & 0.002 & 0.316 & 8.12e-01 & 92.04 \\
 & WSINDy & 0.313 & 0.655 & 6.57e-02 & 20.64 \\
 & PySR & 0.006 & 0.252 & 8.33e-01 & 94.10 \\
 & UNIDE & \textbf{0.941} & \textbf{0.980} & \textbf{3.50e-02} & \textbf{5.61} \\
\midrule
40\% & SINDy & 0.001 & 0.319 & 8.37e-01 & 95.12 \\
 & WSINDy & 0.200 & 0.584 & 8.81e-02 & 27.53 \\
 & PySR & 0.008 & 0.260 & 8.49e-01 & 96.23 \\
 & UNIDE & \textbf{0.938} & \textbf{0.975} & \textbf{3.71e-02} & \textbf{6.86} \\
\midrule
50\% & SINDy & 0.000 & 0.310 & 8.50e-01 & 96.79 \\
 & WSINDy & 0.116 & 0.518 & 1.12e-01 & 33.64 \\
 & PySR & 0.008 & 0.265 & 8.59e-01 & 97.56 \\
 & UNIDE & \textbf{0.864} & \textbf{0.932} & \textbf{4.55e-02} & \textbf{11.54} \\
\midrule
60\% & SINDy & 0.000 & 0.303 & 8.56e-01 & 97.70 \\
 & WSINDy & 0.021 & 0.464 & 1.34e-01 & 34.36 \\
 & PySR & 0.011 & 0.260 & 8.64e-01 & 98.43 \\
 & UNIDE & \textbf{0.911} & \textbf{0.964} & \textbf{4.68e-02} & \textbf{8.91} \\
\midrule
70\% & SINDy & 0.000 & 0.304 & 8.61e-01 & 98.13 \\
 & WSINDy & 0.042 & 0.441 & 1.61e-01 & 42.16 \\
 & PySR & 0.012 & 0.260 & 8.68e-01 & 98.92 \\
 & UNIDE & \textbf{0.810} & \textbf{0.918} & \textbf{7.21e-02} & \textbf{16.17} \\
\midrule
80\% & SINDy & 0.000 & 0.302 & 8.65e-01 & 98.58 \\
 & WSINDy & 0.003 & 0.413 & 1.83e-01 & 43.01 \\
 & PySR & 0.012 & 0.265 & 8.69e-01 & 99.22 \\
 & UNIDE & \textbf{0.924} & \textbf{0.980} & \textbf{5.04e-02} & \textbf{8.04} \\
\bottomrule
\end{tabular}
\end{table}

% let the long noise-sweep table float past this heading (placeins would otherwise leave a half-empty page)
\begingroup\let\FloatBarrier\relax
\section{Generalization Across Noise Levels}
\label{app:perlevel}
\endgroup

Table~\ref{tab:perlevel} evaluates the neural weights trained with noise up to $10\%$ across the noise sweep, alongside matched specialists. The weights are frozen, but the support threshold and pruning tolerance are recalibrated using labeled validation trajectories at each target noise level. This is calibrated weight transfer, not fully frozen zero-shot inference. The support thresholds are $0.7$ at $1\%$ and $5\%$ noise and $0.3$ at $10\%$ and above. The pruning tolerance is $0.05$ through $10\%$, $0.008$ at $20\%$, and zero at $30$--$80\%$. A zero tolerance disables pruning in this implementation. The two configurations coincide at $10\%$ by construction.

\begin{table}[!htbp]
\centering
\small
\setlength{\tabcolsep}{5pt}
\caption{Noise-level generalization: a single model trained with noise up to $10\%$ with frozen neural weights and target-noise validation calibration at every level, versus the models trained at each test level (the main-text configuration). Full 900-trajectory test set, same seed.}
\label{tab:perlevel}
\begin{tabular}{llrrrr}
\toprule
Noise & Training & Exact $\uparrow$ & Term F1 $\uparrow$ & Coef. RMSE $\downarrow$ & Coef. rel\% $\downarrow$ \\
\midrule
1\% & Noise-matched & 0.866 & 0.938 & 5.15e-02 & 10.32 \\
1\% & Trained at 10\% & 0.971 & 0.990 & 2.93e-02 & 3.50 \\
5\% & Noise-matched & 0.877 & 0.942 & 4.27e-02 & 10.18 \\
5\% & Trained at 10\% & 0.952 & 0.981 & 3.23e-02 & 5.12 \\
10\% & Noise-matched & 0.960 & 0.984 & 2.92e-02 & 4.99 \\
10\% & Trained at 10\% & 0.960 & 0.984 & 2.92e-02 & 4.99 \\
20\% & Noise-matched & 0.971 & 0.989 & 3.07e-02 & 4.46 \\
20\% & Trained at 10\% & 0.747 & 0.849 & 5.26e-02 & 21.81 \\
30\% & Noise-matched & 0.941 & 0.980 & 3.50e-02 & 5.61 \\
30\% & Trained at 10\% & 0.667 & 0.779 & 6.60e-02 & 30.85 \\
40\% & Noise-matched & 0.938 & 0.975 & 3.71e-02 & 6.86 \\
40\% & Trained at 10\% & 0.667 & 0.778 & 7.56e-02 & 33.80 \\
50\% & Noise-matched & 0.864 & 0.932 & 4.55e-02 & 11.54 \\
50\% & Trained at 10\% & 0.659 & 0.771 & 1.10e-01 & 37.22 \\
60\% & Noise-matched & 0.911 & 0.964 & 4.68e-02 & 8.91 \\
60\% & Trained at 10\% & 0.652 & 0.766 & 1.43e-01 & 40.32 \\
70\% & Noise-matched & 0.810 & 0.918 & 7.21e-02 & 16.17 \\
70\% & Trained at 10\% & 0.641 & 0.758 & 1.83e-01 & 43.66 \\
80\% & Noise-matched & 0.924 & 0.980 & 5.04e-02 & 8.04 \\
80\% & Trained at 10\% & 0.621 & 0.749 & 2.48e-01 & 47.43 \\
\bottomrule
\end{tabular}
\end{table}

Generalization across noise levels is strongly asymmetric between structure and coefficients. Term selection transfers far beyond the training range: at $80\%$ noise---eight times its training level---the $10\%$-trained model still reaches Term F1 $0.749$, above the three baselines (Figure~\ref{fig:main-results}); toward lower noise it transfers with a gain: at $1\%$ noise it is more accurate than the $1\%$-trained model on every metric (Term F1 $0.990$ versus $0.938$, coefficient error $3.50\%$ versus $10.32\%$), mainly because the $1\%$-trained encoder substitutes $u_x$ for $u\,u_x$ on Burgers. Coefficient accuracy degrades when transferring to higher noise: relative error grows from $5.0\%$ at the training level to $47.4\%$ at $80\%$. The matched-noise model instead has $8.0\%$ relative error at that level. These results motivate the noise-matched protocol in the main text. We apply the protocol uniformly, including at $1\%$, rather than choosing the model per level from test results.

\paragraph{Fully frozen transfer.}
Table~\ref{tab:fully-frozen} holds the denoiser and encoder weights, support threshold $\tau=0.3$, and pruning tolerance $\delta=0.05$ fixed at their $10\%$ source settings. No target-noise validation labels are used to select these settings. Exact match falls from $0.960$ at $10\%$ noise to $0.489$ at $80\%$ noise; mean relative coefficient error rises from $4.99\%$ to $52.00\%$. These results distinguish fully frozen transfer from the target-noise validation calibration in Table~\ref{tab:perlevel} and from matched-noise specialists.

\begin{table}[!htbp]
\centering\small
\caption{Fully frozen transfer of the $10\%$-trained pipeline. Both neural weights, support threshold $\tau=0.3$, and pruning tolerance $\delta=0.05$ remain fixed at source settings. Each row uses all 900 trajectories of the random-trajectory test set; no target-noise validation labels select the settings. The $10\%$ row is the source-level evaluation.}
\label{tab:fully-frozen}
\begin{tabular}{lrrrr}\toprule
Test noise & Exact $\uparrow$ & Term F1 $\uparrow$ & Coef. RMSE $\downarrow$ & Coef. rel\% $\downarrow$ \\ \midrule
1\% & 0.992 & 0.998 & 2.58e-02 & 2.11 \\
5\% & 0.993 & 0.998 & 2.59e-02 & 2.53 \\
10\% & 0.960 & 0.984 & 2.92e-02 & 4.99 \\
20\% & 0.712 & 0.843 & 5.48e-02 & 22.97 \\
30\% & 0.597 & 0.765 & 7.06e-02 & 32.80 \\
40\% & 0.564 & 0.756 & 8.61e-02 & 36.52 \\
50\% & 0.534 & 0.743 & 1.27e-01 & 40.55 \\
60\% & 0.513 & 0.730 & 1.75e-01 & 44.23 \\
70\% & 0.501 & 0.714 & 2.32e-01 & 47.96 \\
80\% & 0.489 & 0.696 & 3.18e-01 & 52.00 \\
\bottomrule\end{tabular}
\end{table}

\section{Per-Family Coefficient-Extrapolation Results at 10\% Observation Noise}
	
Table~\ref{tab:app_per_pde_noise10} breaks the noisy coefficient-extrapolation result of the main text down by PDE family: held-out large-coefficient regimes evaluated under 10\% observation noise, with the best method per metric in bold. With the pruned solve UNIDE reaches exact match of at least 0.91 in five of the six families, including the dispersive and chaotic families where WSINDy fails outright, but only 0.53 on Burgers, where the encoder substitutes $u_x$ for $u\,u_x$ and WSINDy recovers every trajectory; Burgers also accounts for most of UNIDE's relative coefficient error in this setting.

\begin{table}[!htbp]
\centering
\scriptsize
\setlength{\tabcolsep}{4pt}
\caption{Per-family results under coefficient extrapolation at 10\% observation noise. Higher is better for Exact and Term F1; lower is better for coefficient RMSE and relative error. Best within each family block is shown in bold.}
\label{tab:app_per_pde_noise10}
\begin{tabular}{llcccc}
\toprule
PDE & Method & Exact $\uparrow$ & Term F1 $\uparrow$ & Coef. RMSE $\downarrow$ & Coef. rel.\% $\downarrow$ \\
\midrule
Advection & SINDy & 0.000 & 0.172 & 1.77e-01 & 18.54 \\
 & WSINDy & 0.000 & 0.330 & \textbf{3.18e-02} & \textbf{3.34} \\
 & PySR & 0.117 & 0.706 & 2.32e-01 & 24.26 \\
 & UNIDE & \textbf{0.917} & \textbf{0.917} & 8.20e-02 & 8.57 \\
\midrule
Adv.--Diff. & SINDy & 0.000 & 0.304 & 1.26e-01 & 70.35 \\
 & WSINDy & 0.967 & 0.991 & 2.10e-03 & 7.31 \\
 & PySR & 0.000 & 0.453 & 1.57e-01 & 75.68 \\
 & UNIDE & \textbf{1.000} & \textbf{1.000} & \textbf{7.42e-04} & \textbf{1.14} \\
\midrule
Burgers & SINDy & 0.000 & 0.317 & 2.21e-01 & 84.87 \\
 & WSINDy & \textbf{1.000} & \textbf{1.000} & \textbf{4.62e-03} & \textbf{2.19} \\
 & PySR & 0.000 & 0.364 & 2.53e-01 & 91.10 \\
 & UNIDE & 0.533 & 0.781 & 1.40e-01 & 33.35 \\
\midrule
Diffusion & SINDy & 0.000 & 0.167 & 1.40e-02 & 99.64 \\
 & WSINDy & \textbf{1.000} & \textbf{1.000} & 8.75e-04 & 6.21 \\
 & PySR & 0.000 & 0.000 & 1.41e-02 & 100.00 \\
 & UNIDE & \textbf{1.000} & \textbf{1.000} & \textbf{1.46e-04} & \textbf{1.04} \\
\midrule
KdV & SINDy & 0.000 & 0.428 & 3.38e+00 & 93.76 \\
 & WSINDy & 0.000 & 0.503 & \textbf{7.41e-02} & \textbf{4.84} \\
 & PySR & 0.000 & 0.246 & 3.92e+00 & 98.02 \\
 & UNIDE & \textbf{0.983} & \textbf{0.992} & 1.69e-01 & 8.26 \\
\midrule
KS & SINDy & 0.000 & 0.470 & 1.94e+00 & 100.09 \\
 & WSINDy & 0.000 & 0.475 & 1.49e-01 & \textbf{5.99} \\
 & PySR & 0.000 & 0.492 & 2.03e+00 & 98.35 \\
 & UNIDE & \textbf{1.000} & \textbf{1.000} & \textbf{1.35e-01} & 6.94 \\
\bottomrule
\end{tabular}
\end{table}

\section{Per-Family Noise Robustness: Additional Metrics}

\begin{figure}[!htbp]
\centering
\includegraphics[width=\linewidth]{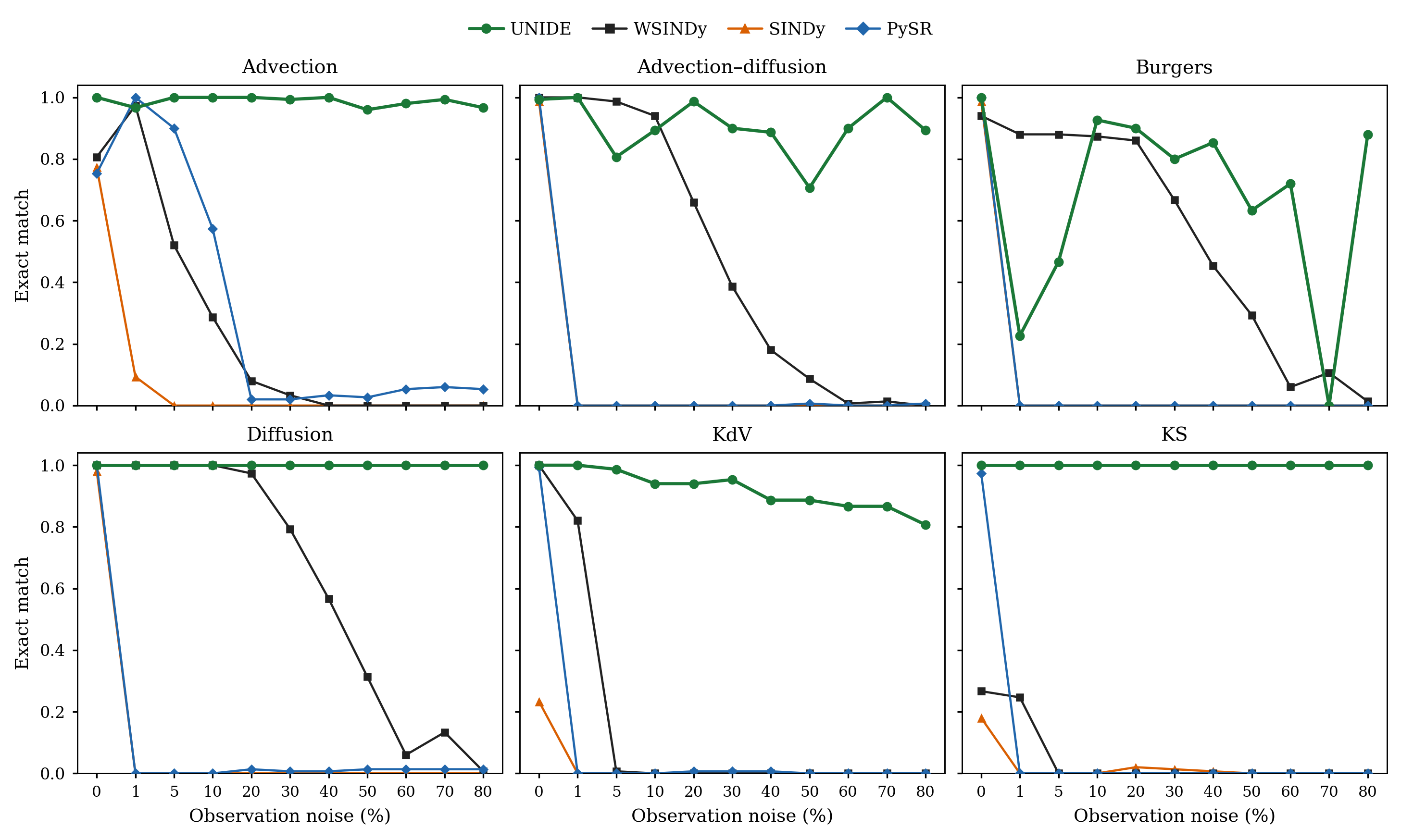}
\caption{Exact match versus observation-noise level for each PDE family. UNIDE keeps near-perfect exact match on advection, diffusion, and KS across the sweep, holds KdV between 0.81 and 1.0 and advection--diffusion between 0.71 and 1.0; on Burgers, UNIDE's exact match is 0.23 at 1\%, 0.47 at 5\%, and 0.00 at 70\%, where the encoder substitutes $u_x$ for $u\,u_x$, and 0.63--0.93 at the other noisy levels.}
\label{fig:per-pde-exact}
\end{figure}

\begin{figure}[!htbp]
\centering
\includegraphics[width=\linewidth]{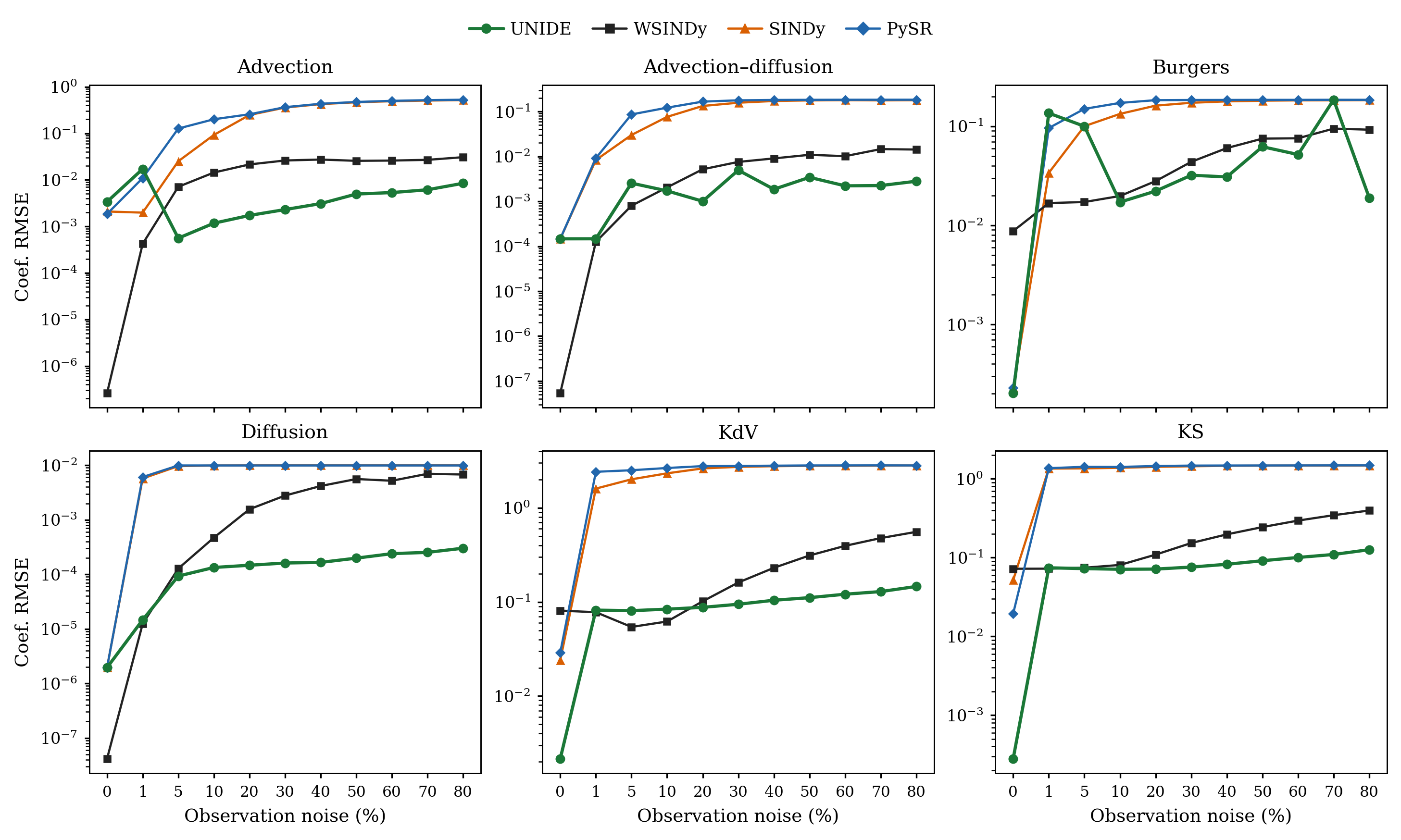}
\caption{Coefficient RMSE (log scale) versus observation-noise level for each PDE family. UNIDE has the lowest coefficient RMSE at seven or more of the ten noisy levels in every family; the pointwise baselines SINDy and PySR inflate sharply once derivatives are taken from noisy data.}
\label{fig:per-pde-rmse}
\end{figure}

\begin{figure}[!htbp]
\centering
\includegraphics[width=\linewidth]{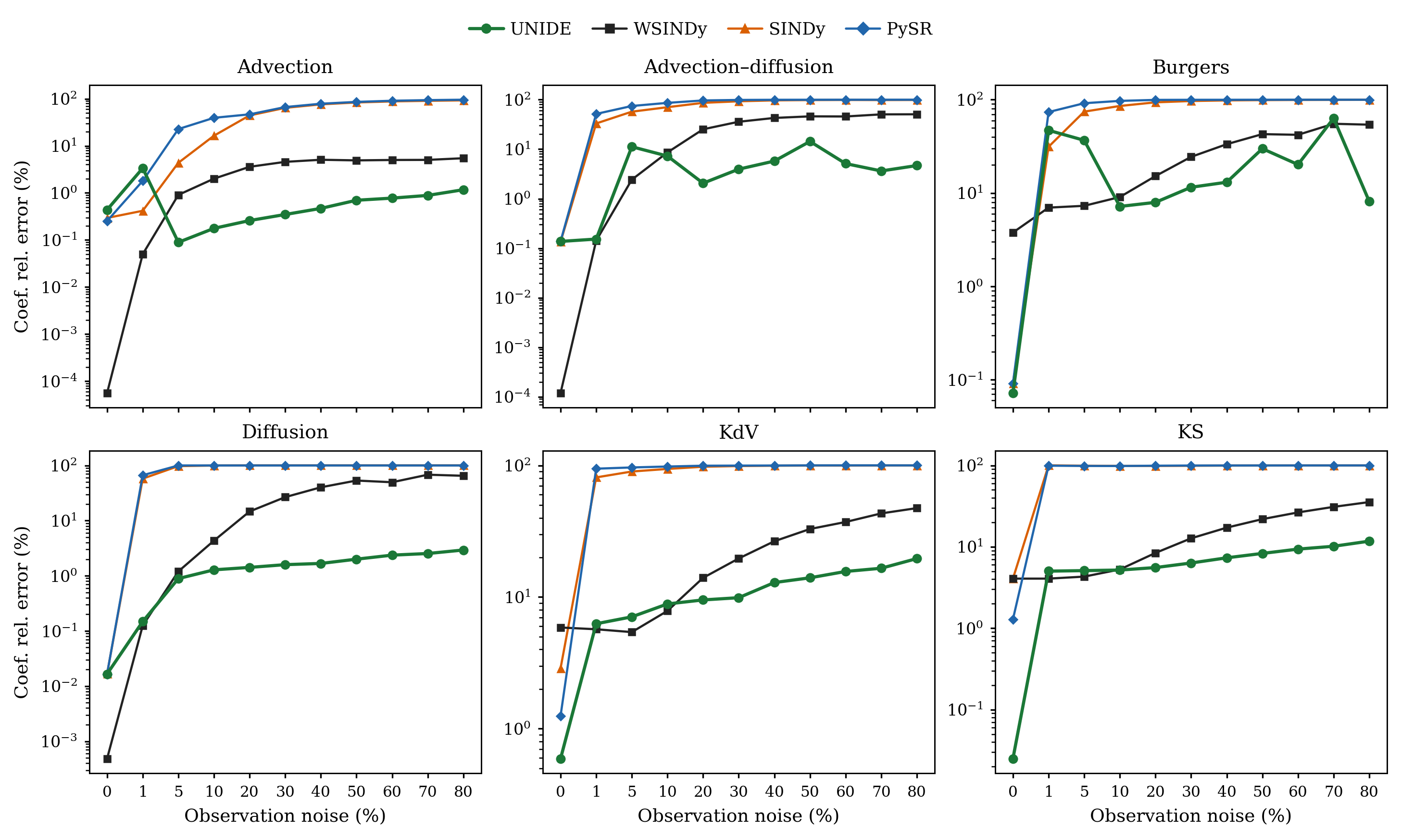}
\caption{Coefficient relative error (\%) versus observation-noise level for each PDE family. UNIDE has the lowest relative error at seven or more of the ten noisy levels in every family, while SINDy and PySR rise toward $100\%$ under noise.}
\label{fig:per-pde-relerr}
\end{figure}

Figure~\ref{fig:per-pde} in the main text reports Term F1 versus observation noise for each PDE family. Figures~\ref{fig:per-pde-exact}--\ref{fig:per-pde-relerr} give the same breakdown for the remaining three metrics --- exact match, coefficient RMSE, and coefficient relative error --- on the full 900-trajectory test set. Methods and colors match Figure~\ref{fig:per-pde}; coefficient metrics are evaluated on the true active terms.

\section{Robustness to Random Seeds}
\label{app:seeds}

The main results use encoder seed 0, the training default for every run, fixed before this study. To measure conditional encoder-seed variability, we retrain with five seeds while keeping the denoiser and common clean warm start fixed, then evaluate all 900 test trajectories at three noise levels. The pruning tolerance is selected on validation separately for each run. Table~\ref{tab:seed-variance} reports means and sample standard deviations. Variability across noise draws is reported separately below.

\begin{table}[!htbp]
\centering
\small
\setlength{\tabcolsep}{6pt}
\caption{Seed variance of UNIDE (mean $\pm$ std over five encoder seeds; per-level denoiser fixed), full 900-trajectory test set.}
\label{tab:seed-variance}
\begin{tabular}{lcccc}
\toprule
Noise & Exact $\uparrow$ & Term F1 $\uparrow$ & Coef. RMSE $\downarrow$ & Coef. rel\% $\downarrow$ \\
\midrule
10\% & $0.926 \pm 0.037$ & $0.966 \pm 0.024$ & $(3.5\pm0.4)\times10^{-2}$ & $7.3 \pm 2.7$ \\
50\% & $0.904 \pm 0.025$ & $0.958 \pm 0.017$ & $(4.2\pm0.3)\times10^{-2}$ & $8.6 \pm 1.9$ \\
80\% & $0.888 \pm 0.052$ & $0.957 \pm 0.034$ & $(5.8\pm1.1)\times10^{-2}$ & $10.8 \pm 4.1$ \\
\bottomrule
\end{tabular}
\end{table}

Term F1 is stable, with a standard deviation of at most $0.034$ across all three levels; the advantage in term recovery over WSINDy is therefore robust to the seed, at matched training and test noise (at $80\%$, F1 $0.957\pm0.034$ versus $0.413$ for WSINDy). Exact match keeps a mean of at least $0.88$ at every level, with a spread of $\pm0.025$ to $\pm0.052$, as expected for the strictest metric. The validation-selected pruning tolerance is identical across seeds at $10\%$ and $50\%$ and spans two grid steps at $80\%$ ($0.004$--$0.012$). Coefficient relative error varies by $\pm1.9$ to $\pm4.1$ points across seeds; its mean stays far below the baselines' at $50\%$ and $80\%$, but at $10\%$ it exceeds WSINDy's ($7.3\%$ versus $6.2\%$). Burgers dominates the spread: its exact match ranges over $0.59$--$0.93$, $0.63$--$0.84$, and $0.00$--$0.88$ across seeds at $10\%$, $50\%$, and $80\%$, the widest range of any family at every level (tied with advection--diffusion at $50\%$). Seed 0 is the strongest of the five at $10\%$ and $80\%$ (e.g., $80\%$ exact match $0.924$ vs.\ mean $0.888\pm0.052$) and the weakest at $50\%$ ($0.864$ vs.\ $0.904\pm0.025$), so both the seed-0 values and the five-seed summaries are needed to describe these experiments.

\paragraph{Observation-noise replication.}
\label{app:noise-replication}
Table~\ref{tab:noise-replication} evaluates three observation-noise realizations at each of $10\%$ and $80\%$ noise on the same 900 trajectories of the random-trajectory test set. Exact match ranges from $0.954$ to $0.961$ at $10\%$ noise and from $0.921$ to $0.926$ at $80\%$. Only the noise draw changes: the trained weights and support threshold are fixed, and the pruning tolerance is reselected on validation (the selected tolerances at $10\%$ noise are $0.05$, $0.02$, and $0.05$).

\begin{table}[!htbp]
\centering\small
\caption{Observation-noise replication on the same 900 trajectories of the random-trajectory test set per row. The denoiser and encoder weights and source support threshold ($0.3$ at $10\%$ noise, $0.7$ at $80\%$) are fixed. Pruning tolerance is reselected on validation for each realization; the selected values are $0.05$, $0.02$, and $0.05$ at $10\%$ and $0.012$ for all three $80\%$ runs. A seed offset plus the sample index determines the noise draw; offset 1000 reproduces the original realization.}
\label{tab:noise-replication}
\begin{tabular}{lrrrrr}\toprule
Noise & Seed offset & Exact $\uparrow$ & Term F1 $\uparrow$ & Coef. RMSE $\downarrow$ & Coef. rel\% $\downarrow$ \\ \midrule
10\% & 1000 & 0.960 & 0.984 & 2.92e-02 & 4.99 \\
10\% & 2000 & 0.954 & 0.974 & 3.20e-02 & 5.95 \\
10\% & 3000 & 0.961 & 0.981 & 3.06e-02 & 5.18 \\
80\% & 1000 & 0.924 & 0.980 & 5.04e-02 & 8.04 \\
80\% & 2000 & 0.921 & 0.981 & 4.95e-02 & 8.09 \\
80\% & 3000 & 0.926 & 0.982 & 4.74e-02 & 7.91 \\
\bottomrule\end{tabular}
\end{table}

\end{document}